\documentclass[12pt]{article}

\usepackage{color}

\usepackage[utf8]{inputenc}
\usepackage{epsfig}
\usepackage{amsmath}
\usepackage{amsfonts}
\usepackage{amsthm}
\usepackage{blkarray}
\usepackage{mathrsfs}
 \usepackage{stmaryrd}
\usepackage{graphicx} 
\usepackage{subfigure}
\usepackage{empheq}
 \usepackage{indentfirst}
\usepackage{mathenv}
\usepackage{fancybox}
\usepackage[T1]{fontenc}
\usepackage{lmodern}
\usepackage[toc,page]{appendix}
\usepackage[british,UKenglish,USenglish,english,american]{babel}
\usepackage{pifont}
\usepackage{dsfont}

\usepackage{float}
\usepackage{placeins}

        \theoremstyle{plain}

        \theoremstyle{remark}
        \newtheorem{remark}{\bf Remark}[section]
        \theoremstyle{remark}
        
        \theoremstyle{remark}

\newcommand{\demi}{\frac{1}{2}}

\newcommand{\R}{{\mathbb{R}}}

\newcommand{\dt}{\partial_t}
\newcommand{\dx}{\partial_x}

\newcommand{\x}{{\bf x}}

\newcommand{\vit}{{\bf v}}
\newcommand{\um}{{\bf u}}

\newcommand{\Maxw}{M}

\newcommand{\ximdemi}{x_{i-\frac{1}{2}}}
\newcommand{\xipdemi}{x_{i+\frac{1}{2}}}

\newcommand{\Fipdemin}{{\cal F}_{i+\frac{1}{2},k}^n}

\newcommand{\fipdeminp}{{f}_{i+\frac{1}{2},k}^{n,+}}
\newcommand{\fipdeminm}{{f}_{i+\frac{1}{2},k}^{n,-}}

\newcommand{\modif}[1]{{\color{blue}{#1}}}

\newcommand{\yee}{O2-flux}
\newcommand{\slopenolim}{O2-slope-nolim}
\newcommand{\slope}{O2-slope}
\newcommand{\slopeCLzero}{O2-slope-BC-O1}
\newcommand{\DG}{DG}

\newcommand{\first}{O1}
\newcommand{\NS}{NS}

\begin{document}

\begin{center}

{\bf Numerical boundary conditions in Finite Volume and Discontinuous Galerkin schemes for
the simulation of rarefied flows along solid boundaries}

\vspace{1cm}
C. Baranger$^1$, N. H\'erouard$^{1,2}$, J. Mathiaud$^{1,2}$, L. Mieussens$^2$

\bigskip
$^1$CEA-CESTA\\
15 avenue des Sabli\`eres - CS 60001\\
33116 Le Barp Cedex, France\\
{ \tt(celine.baranger@cea.fr, julien.mathiaud@cea.fr)}\\

\bigskip
$^2$ Bordeaux INP, Univ. Bordeaux, CNRS, INRIA, IMB, UMR 5251, F-33400 Talence, France.\\

{ \tt(Nicolas.Herouard@math.u-bordeaux.fr, Luc.Mieussens@math.u-bordeaux.fr)}

\end{center}

\paragraph{Abstract:} We present a numerical comparison between two
standard finite volume schemes and a discontinuous Galerkin method
applied to the BGK equation of rarefied gas dynamics. We pay a
particular attention to the numerical boundary conditions in order to
preserve the rate of convergence of the method. Most of our analysis
relies on a 1D problem (Couette flow), but we also
present some results for a 2D aerodynamical flow.

\bigskip


\tableofcontents

\section{Introduction}

The kinetic theory of Rarefied gas flows describes the behavior of a
gas in a system which length is of same order of magnitude as the mean
free path of the gas molecules. Numerical simulations of these flows
are fundamental tools, for instance they are used in aerodynamics to
estimate heat fluxes at the surface of a re-entry space vehicle at
high altitudes or to estimate the attenuation of a micro-accelerometer
by the surrounding gas in a micro-electro-mechanical system. One last
example is when one wants to estimate the pumping speed or compression
rate of a turbo-molecular pump.

For such problems, it is important to take into account the
interactions of the gas with the solid boundaries of the system. They
are modeled by boundary conditions like the diffuse or specular
reflections. In this work, our first aim is to improve the accuracy
of our in-house kinetic solver, which is based on a finite volume
approach (FV), in particular for the computation of parietal fluxes. Finite
volume schemes are often used in rarefied gas dynamics (RGD), see for
instance the work of~\cite{T_cicp, luc_jcp, Baranger,DM_2016, XH2010}. Our scheme uses a standard extension to
second order accuracy by using the non linear flux limiter approach of
Yee~\cite{yee}. We observe that the accuracy of this scheme decreases
to first order at solid boundaries, because the boundary conditions
are not discretized properly. The first goal of this paper is to
analyze this problem, both analytically and numerically. We show that
our flux limiter can hardly be compatible with second order
discretization of reflection boundary conditions.

Then we study another kind of second order extension based on linear
reconstruction and nonlinear slope limiters~\cite{leveque}. In this
case, we are able to propose a discretization of the boundary
conditions (inspired by~\cite{Filbet201343,guigui}), based on extrapolations that are consistent with the slope
limiter technique, which preserves second order accuracy up to solid
boundaries. On numerical tests, the improvement of the accuracy of
this method is spectacular, for 1D and 2D problems.

Finally, we want to compare this modified finite volume scheme to
the popular Discontinuous Galerkin method (DG). Incidentally, this
method has originally been made for neutron transport (hence a kinetic
equation, like the Boltzmann equation of RGD) in the 1970's by Reed
and Hill~\cite{Reed_Hill_1973}. Such schemes have many interesting
properties, and they became more and more popular in the past decades
in many different fields. In particular, the DG method has been used
for charged particles~\cite{gamba}, and more recently for the BGK
equation~\cite{SAC_2012}. According to this last paper, the DG
method would be more efficient that the FV method, at least for the
test cases used by the authors. In this paper, we show that the
discretization of reflection conditions with a DG method is very
simple. We compare a simple second order DG scheme to our two FV
schemes, mainly to investigate their accuracy at solid
boundaries. This study shows the order of accuracy of both DG and FV
method (with slope limiters) are not decreased at solid boundaries.

The outline of this article is the following. In
section~\ref{sec:BGK}, we present the basic BGK equation of RGD. 
In~\ref{sec:numsche}, we analyze the two versions of our FV scheme and
a simple DG scheme. These schemes are compared for several test cases
in section~\ref{sec:Num}.

\section{The BGK equation}
\label{sec:BGK}

In kinetic theory, a monoatomic gas is described by a particle distribution
function $f(t,\x,\vit)$ defined such as $f(t,\x,\vit)\textrm{d}\x \textrm{d}\vit$ is the mass
of molecules that at time $t$ are located in an elementary space volume
$d\x$ centered in $\x=(x,y,z)$ and have a velocity in an elementary
volume $d\vit$ centered in $\vit=(v_x,v_y,v_z)$.\
The macroscopic quantities associated to $f$ such  as mass density $\rho$,
momentum $\rho \um$ and total energy $E$ are defined as the five first
moments of $f$ with respect to the velocity variable, namely:
\begin{equation}  \label{eq-mts}
(\rho(t,\x),\rho\um(t,\x),E(t,\x))=\int_{\R^3}(1,\vit,\demi|\vit|^2)f(t,\x,\vit)
\, \textrm{d}\vit.
\end{equation}
The temperature $T$ of the gas and its pressure $p= \rho R T$ are defined by the relation
\begin{equation}  \label{eq-temp}
E=\demi\rho|\um|^2+\frac{3}{2}\rho R T
\end{equation}
where $R$ is the gas constant
defined as the ratio between the Boltzmann constant and the molecular
mass of the gas.
When the gas is in a thermodynamical equilibrium state, it is well known that
the distribution function $f$ is a Gaussian function
$\Maxw[\rho,\um,T]$ of $\vit$, called
Maxwellian distribution, that depends only on 
macroscopic quantities and satisfies relations~(\ref{eq-mts}):
\begin{equation}  \label{eq-maxw}
\Maxw[\rho,\um,T]=\frac{\rho}{(2\pi RT)^\frac{3}{2}}\exp(-\frac{|\vit-\um|^2}{2RT}).
\end{equation}

The evolution of a  gas at thermodynamical nonequilibrium is described by the following Boltzmann equation:
\begin{equation}  \label{eq-kin}
\dt f + \vit \cdot \nabla_\x f = Q(f),
\end{equation}
which means that the total variation of $f$ (described by the
  left-hand side) is due to  collisions between molecules
  ($Q(f)$ is the collision kernel). The most realistic collision
kernel is the Boltzmann operator but its use is still very computationally
expensive. Here we choose the simpler BGK
model~\cite{bgk,welander}
\begin{equation}  \label{eq-bgk}
Q(f) = \frac{1}{\tau}(\Maxw[\rho,\um,T]-f)
\end{equation}
which makes $f$ relax
towards the local equilibrium corresponding to the macroscopic
quantities defined by~(\ref{eq-mts}).
The relaxation time  defined
as $\displaystyle{\tau=\frac{\mu}{\rho RT}}$ ($\mu$ being the viscosity of the
gas)  is chosen to recover the correct viscosity in the Chapman-Enskog~expansion
 for Navier-Stokes equations. 


The interactions of the gas with solid boundaries are described with
\textcolor{black}{a diffuse reflection model.}
If the boundary has a temperature $T_w$ and a 
velocity $\um_w$ so that $\um_w \cdot {\bf n}(\x) = 0$, where ${\bf n}(\x)$ is the normal to the wall at point $\x$ directed into the gas, a molecule that
collides with the boundary is  re-emitted   with a random velocity normally
distributed around $\um_w$ ($RT_w$ being the  variance of the normal distribution). This reads
\begin{equation}\label{eq-BCdiff} 
  f(t,\x,\vit)_{|\Gamma} =  \sigma_w(f)\Maxw[1,\um_w,T_w](\vit)\ \  \textrm{if } \vit\cdot {\bf n}(\x) > 0.
\end{equation}

The parameter $\sigma_w(f)$ \textcolor{black}{is the mass density of
  reflected particles, and it}
is defined so that there is no normal mass flux across the boundary (all
the molecules are re-emitted), \textcolor{black}{see~\cite{cercignani}}.
Namely, that is,  
\begin{equation}  \label{eq-sigmadiff}
\sigma_w(f)=-\left(\int_{\vit\cdot {\bf
    n}(\x) < 0} \vit\cdot {\bf n}(\x) f(t,\x,\vit)  \, \textrm{d}\vit\right)\left(\int_{\vit\cdot {\bf
    n}(\x) > 0} \vit\cdot {\bf n}(\x)  \Maxw[1,\um_w,T_w](\vit) \, \textrm{d}\vit\right)^{-1}.
\end{equation}
There are other reflection models, like the Maxwell model with partial
accommodation, but they are not used in this work.

The numerical resolution of the BGK equation is done by a
deterministic method using a discrete velocity approach (see
\cite{luc_jcp,MieussensL,Mieussens200183}, and the extension to
polyatomic gases~\cite{luc_dubroca}). The velocity variable is
replaced by discrete values $\vit_k$ of a Cartesian
grid (see \cite{Baranger} for the use of locally refined velocity
grids). The continuous distribution $f$ is then replaced by its
approximation at each point $\vit_k$, and we get the following
discrete velocity BGK equation
\begin{equation}
\frac{\partial  f_k}{\partial t} +\vit_k \cdot \nabla_\x  f_k=
\frac {1}{\tau}(\mathcal{M}_k[\rho,\um,T]  - f_k), \label{BGKdiscret}
\end{equation}
where $f_k(t,\x) \approx f(t,\x,\vit_k)$. The discrete function
$\mathcal{M}_k$ is an approximation of the equilibrium functions
$\displaystyle{\Maxw}$ defined in \eqref{eq-maxw} at point
$\vit_k$. This discrete equilibrium is defined in order to satisfy the
conservation of mass, momentum and energy (see~\cite{luc_jcp}). The
macroscopic quantities are now defined by the quadrature rule
\begin{equation}  \label{eq-mtsk}
(\rho(t,\x),\rho\um(t,\x),E(t,\x))=\sum_k(1,\vit_k,\demi|\vit_k|^2)f_k(t,\x)\, \omega_k,
\end{equation}
where $\omega_k$ are the quadrature weights: \textcolor{black}{ for a
  Cartesian velocity grid of steps $\Delta v_{x}$, $\Delta v_{y}$, $\Delta v_{z}$, the weights are
  given by $\omega_k=\Delta v_x \Delta v_y \Delta v_z$}. Moreover, the
boundary condition~(\ref{eq-BCdiff}) is replaced by
\begin{equation}\label{eq-BCdiffk} 
  f_k(t,\x)_{|\Gamma} =  \sigma_w(f){\cal M}_k[1,\um_w,T_w]\ \  \textrm{if } \vit_k\cdot {\bf n}(\x) > 0,
\end{equation}
\textcolor{black}{where $\sigma_w(f)$ is defined as
  in~(\ref{eq-sigmadiff}) in which the integrals are replaced by
  quadratures: }
\begin{equation}  \label{eq-sigmadiffk}
\sigma_w(f)=-\left(\sum_{\vit_k \cdot {\bf
    n}(\x) < 0}  \vit_k\cdot {\bf n}(\x)  f_k(t,\x)\, \omega_k\right)\left(\sum_{\vit_k\cdot {\bf
    n}(\x) > 0} \vit_k\cdot {\bf n}(\x) {\cal M}_k[1,\um_w,T_w] 
\, \omega_k \right)^{-1}.
\end{equation}
\textcolor{black}{Again this ensures a zero mass flux across the solid wall.}

   \section{Numerical schemes and boundary conditions}
      \label{sec:numsche}
For the sake of simplicity, the numerical scheme will be presented in
the 1D configuration, in the domain $[0,1]$. Equation \eqref{BGKdiscret} is then written
\begin{equation}
\frac{\partial  f_k}{\partial t} +v_k \frac{\partial  f_k}{\partial x}=
\frac {1}{\tau}(\mathcal{M}_k[\rho,\um,T]  - f_k), \label{BGKdiscret1D}
\end{equation}
where $v_k$ is the first component of $\vit_k$. Since $\vit_k\cdot
{\bf n}(\x) = v_k$ (at $x=0$), the boundary
condition~(\ref{eq-BCdiffk})--(\ref{eq-sigmadiffk}) now reads
\begin{equation}\label{eq-BCdiffk_1d} 
  f_k(t,0) =  \sigma_w(f){\cal M}_k[1,\um_w,T_w]\ \  \textrm{if } v_k > 0,
\end{equation}
where
\begin{equation}  \label{eq-sigmadiffk_1d}
\sigma_w(f)=-\left(\sum_{v_k  < 0} v_k  f_k(t,0) \,
  \omega_k\right)\left(\sum_{v_k > 0}v_k {\cal M}_k[1,\um_w,T_w]  \, \omega_k \right)^{-1}.
\end{equation}

For the following, it is important to note that boundary
condition~(\ref{eq-BCdiffk_1d})--(\ref{eq-sigmadiffk_1d}) implies a zero
mass flux $\phi$ at the boundary. Indeed,~(\ref{eq-BCdiffk_1d}) yields: 
\begin{equation*}
\begin{split}
\phi & := \sum_k v_k f_k(t,0) \, \omega_k \\
&  =  \sum_{v_k  < 0} v_k  f_k(t,0) \,  \omega_k  + \sum_{v_k  < 0} v_k
\sigma_w(f){\cal M}_k[1,\um_w,T_w] \,  \omega_k 
 = 0
\end{split}
\end{equation*}
from~(\ref{eq-sigmadiffk_1d}). In particular, this property ensures a
global mass conservation for any internal flow.

\subsection{Finite Volume schemes}
\label{subsec:VF}
Finite Volume schemes approximate the solution of a given problem by
integrating the equation on each cell of a mesh. The integration of
the advection term results in a flux at cell interfaces: the numerical
flux. The accuracy of the scheme depends on the accuracy of this
flux. Moreover, this flux is closely linked with the boundary
conditions. This is why, in this section, we study in detail three different finite
volume schemes and their properties close to a solid boundary. Note
that the collision term $Q(f)$ is local in space and its
discretization is consequently the same for all Finite Volume schemes:
therefore, even if we take the BGK equation as an
example, the study performed in this section can be applied to
every collision operator.

First, we assume we have a mesh of $i_{max}+1$ nodes $x_{i+\demi}$ with steps
$\Delta x_i = x_{i+\demi}-x_{i-\demi}$, for $i=1$ to $i=i_{max}$. The discrete time variable is
$t_n$ with a time step $\Delta t_n$. We first integrate
equation~(\ref{BGKdiscret1D}) in a cell $[x_{i-\demi},x_{i+\demi}]$
between $t_n$ and  $t_{n+1}$ and divide by $\Delta t_n$ to obtain
\begin{equation}\label{eq-fvol} 
  \frac{f^{n+1}_{i,k}-f^{n}_{i,k}}{\Delta t_n} + \frac{1}{\Delta
    x_i}\left( {\cal F}^n_{i+\demi,k} - {\cal F}^n_{i-\demi,k} \right)
  = Q_k(f^n_{i}),
\end{equation}
where $f^{n}_{i,k}$ is an approximation of the average of $f_k(t_n,x)$
in the cell $i$. The numerical flux ${\cal F}^n_{i+\demi,k}$ is an
approximation of the integral $\int_{t^n}^{t^{n+1}}v_k
f_k(t,x_{i+\demi})\, dt$, and can be interpreted as a flux across the
cell interface $x_{i+\demi}$ between cells $i$ and $i+1$.

In the literature, there are many different constructions of this
numerical flux, and each one leads to a different scheme. The simplest
one is the upwind scheme (see section~\ref{subsec:VF1} below), which
gives a first order scheme. There are several ways to increase the
order of a this scheme, all of them related to a modification of the
numerical flux. One solution is to reconstruct the solution inside the
cells of the mesh, in order to define new values on the interfaces. Another one
is to limit a high order numerical flux at the cell interfaces. Here,
we study one method of each category, and in particular their
properties at the solid boundary. We first present the standard first
order upwind scheme. Then we study a second order flux limiter method
and show it decreases to first order at the boundary. Then we present
the linear reconstruction method, which gives a correct second-order
approximation of the BGK problem up to the boundary.

\modif{Finally, note that, for simplicity, all our schemes are presented with a
  standard forward (explicit) Euler time discretization. In
  practice, since we are interested in steady flows only, our numerical tests are made with a linearized backward
  (implicit) Euler method, which is standard in aerodynamics
  (see~\cite{luc_jcp,Baranger}). However, the analysis made below does not depend
  on the time discretization.}

    \subsubsection{A first order scheme finite volume scheme}
\label{subsec:VF1}

The simplest numerical flux is the upwind flux:
\begin{equation*}
\Fipdemin =v_k^{+}{f}_{i,k}^n+v_k^{-}{{f}_{i+1,k}^n},
\end{equation*}
where $v_k^{\pm} = (v\pm|v|)/2$ denotes the positive and negative
parts of $v_k$. It is well known that the corresponding scheme is
first order accurate in space, with a strong numerical diffusion on
coarse meshes.

Assume that the first cell of the mesh $[x_{\demi},x_{\frac32}]$ is adjacent to a solid
wall. Then the numerical flux
${\cal F}^n_{\demi,k}$ approximates the flux across the solid boundary. It
requires an artificial value $f^n_{0,k}$, that is generally
interpreted as the value of the distribution function in a {\it ghost
  cell} $[x_{-\demi},x_{\demi}]$ inside the solid wall, adjacent to
the first cell. This
value must be defined so as to (a) account for the boundary
condition~(\ref{eq-BCdiffk_1d})--(\ref{eq-sigmadiffk_1d}), and (b) satisfy a zero
mass flux across the wall, that is to say $\sum_k{\cal  F}^n_{\demi
  ,k} \omega_k=0$. 

These constraints can be satisfied as follows. First, the value of
$f^n_{0,k}$ is defined for outgoing velocities $v_k<0$ by using a
zeroth order extrapolation of the value in $\Omega_1$, namely
\begin{equation*}
  f^n_{0,k} = f^n_{1,k}.
\end{equation*}
Then, the value of $f^n_{0,k}$ for incoming velocities is defined by
using the boundary condition~(\ref{eq-BCdiffk_1d}-\ref{eq-sigmadiffk_1d}),
that is to say:
\begin{equation*}
  f^n_{0,k} = \sigma_w(f^n_{0}){\cal M}_k[1,\um_w,T_w], 
\quad \text{ where } \quad 
\sigma_w(f^n_{0})  = - \left({\sum_{l}}v_l^{-}
       {f}_{0,l}^{n}\omega_l\right)/\left(\sum_{l}v_l^{+}
       \mathcal{M}_l[1,\um_w,T_w]\omega_l\right).
\end{equation*}
It can easily seen that this definition ensures the mass
conservation. Indeed, the mass flux across the wall is 
\modif{
\begin{equation*}
\begin{split}
\sum_k {{\cal F}_{1/2,k}^n}\omega_k 
& = \sum_k \left(v_k^{+}{f}_{0,k}^n+v_k^{-}{f}_{1,k}^n\right)\omega_k
\\
& = \sum_k \left( -v_k^{+}\displaystyle{\frac{{\sum_{l}}v_l^{-} {f}_{0,l}^{n}\omega_l}{\sum_{l}v_l^{+}
      \mathcal{M}_l[1,\um_w,T_w]\omega_l}\mathcal{M}_k[1,\um_w,T_w]}\right) \omega_k 
 + \sum_k \left( v_k^{-} {f}_{1,k}^{n}\right) \omega_k \\
& =-\sum_l\left(v_l^{-} {f}_{0,l}^{n}\right) \omega_l 
 + \sum_k \left( v_k^{-} {f}_{1,k}^{n}\right) \omega_k \\
& = 0,
\end{split}
\end{equation*}
}
since ${f}_{0,k}^{n} = {f}_{1,k}^{n}$ for $v_k<0$.



   \subsubsection{Second order finite volume scheme with flux limiters}
   \label{schemyee}
   
   In this kind of scheme, a nonlinear limitation is applied to a
   centered numerical flux (which is second order accurate) to ensure
   a stability property. Here, we use the Yee limiter, which ensures
   that the scheme is ``Total Variation Diminishing'' (TVD) (see \cite{yee}):
\begin{equation}\label{eq-yee} 
\left\{
      \begin{aligned}
&\Fipdemin = v_k^{+}{{f}_{i,k}^n}+v_k^{-}{{f}_{i+1,k}^n}+|v_k|\displaystyle{\frac{1}{2}}{\Phi_{i+1/2,k}^n},\\
&{\Phi_{i+1/2,k}^n} = {\rm minmod}\displaystyle{\left(\Delta {f}_{i-1/2,k}^n,\Delta {f}_{i+1/2,k}^n,\Delta {f}_{i+3/2,k}^n\right)},\\
&\Delta {f}_{i+1/2,k}^n = {f}_{i+1,k}^n-{f}_{i,k}^n.
      \end{aligned}
\right.
\end{equation}\\
where the minmod function is defined by
\begin{equation*}
\textrm{minmod}(x,y,z) =
\left\{
\begin{aligned}
& \textrm{sgn}(x)\ \textrm{min}(|x|,|y|,|z|)\ \ \textrm{if}\ \textrm{sgn}(x)=\textrm{sgn}(y)=\textrm{sgn}(z),\\
&0 \  \textrm{in other cases}.
\end{aligned}
\right.
\end{equation*}

Like in the first order scheme, some ghost cell values have to be
defined to compute the numerical flux at the wall interface ${\cal
  F}^n_{\demi ,k}$. Indeed, the first order part of the flux requires the value of
$f^n_{0,k}$ for incoming velocities ($v_k>0$), while the limiter needs
$f^n_{0,k}$ and $f^n_{-1,k}$ for every velocities. Again, these ghost cell
values have to be defined so as to account for the boundary
condition~(\ref{eq-BCdiffk_1d})--(\ref{eq-sigmadiffk_1d})
to satisfy a zero
mass flux across the wall.

A simple idea that satisfies these constraints is the following. For
outgoing velocities, we define $f^n_{0,k}$ and $f^n_{-1,k}$ by a
zeroth order extrapolation from the first boundary cell: 
\begin{equation*}
f^n_{-1,k} =   f^n_{0,k}  = f^n_{1,k}.
\end{equation*} 
For incoming velocities, we define $f^n_{0,k}$  like in the first
order scheme, that is to say: 
\begin{equation*}
  f^n_{0,k} = \sigma_w(f^n_{0,k}){\cal M}_k, 
\quad \text{ where } \quad 
\sigma_w(f^n_{0,k}) = - \displaystyle{\left({\sum_{l}}v_l^{-}
       {f}_{0,l}^{n}\omega_l\right)\left(\sum_{l}v_l^{+}
       \mathcal{M}_l\omega_l\right)^{-1}}.
\end{equation*}
And again, $f^n_{-1,k}$ is defined by zeroth order extrapolation with
$f^n_{-1,k} =   f^n_{0,k} $.

This definition implies that the numerical flux at the wall interface ${\cal
  F}^n_{\demi ,k}$ reduces to the first order flux
$v_k^{+}{{f}_{0,k}^n}+v_k^{-}{{f}_{1,k}^n}$, since $\Delta
f^n_{-\demi,k}=f^n_{0,k}-f^n_{-1,k} =0$ by construction, and hence the
limiter $\Phi_{1/2,k}^n$ is zero. Consequently, the mass conservation
is satisfied (see section~\ref{subsec:VF1}).
At the same time, the drawback of this approach is that the numerical
flux at the wall is only first order accurate, while it is second
order inside the computational domain. This lower accuracy at the wall
can be clearly observed in our numerical tests (see
section~\ref{sec:Num}).

To increase the accuracy of the scheme, a natural idea is to use a
first order extrapolation to define $f^n_{0,k}$ and $f^n_{-1,k}$ for
outgoing velocities, and then to use the boundary conditions to define
these values for incoming velocities. This might give a second order
flux, but we have not been able to find a definition for which the
conservation property still holds. In our opinion, it is unlikely that
such a definition exist \textcolor{black}{(see an explanation at the end of section~\ref{reconslin})}.

Consequently, we believe that such a scheme cannot be at the same time
second order up to the boundary and conservative. This is why we
propose to consider another kind of second order scheme in the
following section.

\subsubsection{Second order finite volume second order scheme based on
  linear reconstruction}
\label{reconslin}

In the finite volume schemes using a linear reconstruction, the
numerical fluxes appearing in relation $\eqref{eq-fvol}$ are still
defined as upwind fluxes, but with new values at the cells interfaces,
obtained with a linear reconstruction of the distribution $f_k$ inside
each cell~\cite{Leveque_book}:
\begin{equation}\label{eq-rec} 
\left\{
      \begin{aligned}
&\Fipdemin = v_k^{+}\fipdeminm+v_k^{-}\fipdeminp,\\
&\fipdeminm={{f}_{i,k}^n}+\frac{\Delta x_{i}}{2} {{\delta}_{i,k}^n},\\
&\fipdeminp={{f}_{i+1,k}^n}-\frac{\Delta x_{i+1}}{2} {{\delta}_{i+1,k}^n},\\
      \end{aligned}
\right.
\end{equation}
where ${{\delta}_{i,k}^n}$ is the slope of the linear reconstruction
of $f$ in cell $i$. The most classical way to compute this slope is a
least square method, which gives, for a regular mesh, the centered slope
\begin{equation}
{\delta}_{i,k}^n = \displaystyle{\frac{{{f}_{i+1,k}^n}-{{f}_{i-1,k}^n}}{2\Delta x_{i}}}. \label{moindrecarre}
\end{equation}
To obtain a TVD scheme, the slope must be limited to avoid the
creation of new extrema. Here, we use the MC slope limiter defined on each cell by
\begin{equation}
 {{\delta}_{i,k}^{n,lim}} = {\rm minmod}\displaystyle{\left({\delta}_{i,k}^n,2 \alpha_i \frac{{{f}_{i,k}^n}-{{f}_{i-1,k}^n}}{\Delta x_{i}},2 \alpha_i \frac{{{f}_{i+1,k}^n}-{{f}_{i,k}^n}}{\Delta x_{i}}\right)}, \label{limgeneral}
\end{equation}
where $\alpha_i$ is a free parameter between 0 and 1~\cite{Leveque_book}.

When this scheme is applied to the first cell $i=1$, the numerical
flux ${\cal F}^n_{\demi,k}$ is 
\begin{equation*}
  {\cal F}^n_{\demi,k}= v_k^{+}f^{n,-}_{\demi,k} + v_k^{-}f^{n,+}_{\demi,k}.
\end{equation*}
It uses two values of the distribution at the solid interface:
$f^{n,+}_{\demi,k}$ is the value on the right side of the wall,
required for outgoing velocities $v_k<0$, given
by the linear reconstruction
\begin{equation}\label{eq-fdemiplus} 
  f^{n,+}_{\demi,k} = {{f}_{1,k}^n}-\frac{\Delta x_{1}}{2} {{\delta}_{1,k}^{n,lim}},
\end{equation}
while $f^{n,-}_{\demi,k}$ is the value on the left side of the wall,
required for incoming velocities $v_k>0$, given
by the linear reconstruction
\begin{equation}\label{eq-fdemimoins} 
  f^{n,-}_{\demi,k} = {{f}_{0,k}^n}+\frac{\Delta x_{0}}{2} {{\delta}_{0,k}^{n,lim}},
\end{equation}
where the slopes are defined
by~(\ref{moindrecarre})--(\ref{limgeneral}) with $i=1$ and $i=0$,
respectively. This
requires to define two ghost cell values $f^n_{0,k}$ and $f^n_{-1,k}$, and the
corresponding cell size $\Delta x_0$. Again, they have to be defined so
as to account for the boundary condition and to preserve the
conservation property.

First, note that for the conservation property holds, it is sufficient that
the value of $f$ at the wall interface satisfies 
\begin{equation}\label{eq-fdemicl} 
  f^{n,-}_{\demi,k} = \sigma_w(f^{n,+}_{\demi,k}) {\cal M}_k.
\end{equation}
Consequently, we have to construct the ghost cell values so
that relations~(\ref{eq-fdemiplus})--(\ref{eq-fdemicl}) hold. This
can be done as follows.

\begin{figure}[p]
\centering
\includegraphics[scale=0.5]{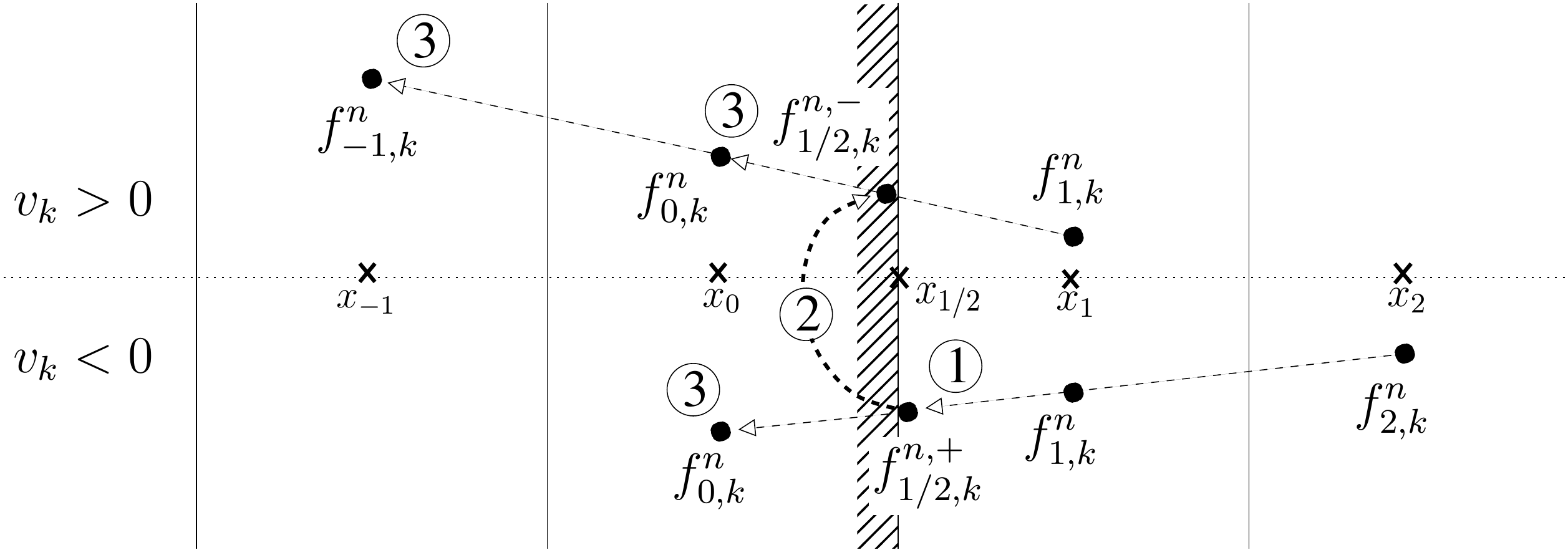}
\caption{Finite volume scheme with a slope limiter: definition of the ghost cell values
  $f^n_{0,k}$ and $f^n_{-1,k}$ by using linear extrapolations and the
  boundary condition.}
\label{O2}
\end{figure}

The method, inspired by ~\cite{Filbet201343,guigui}, consists in first
defining the wall interface values $f^{n,\pm}_{\demi,k} $ by
extrapolation (for outgoing velocities) and by the boundary condition
(for incoming velocities), and then in defining the ghost cell values
$f^n_{0,k}$ and $f^n_{-1,k}$ so that
~(\ref{eq-fdemiplus})--(\ref{eq-fdemimoins}) hold. This method is
detailed below in three steps, and summarized in
figure~\ref{O2}. Note that in the following, we will assume that
the ghost cells have the same size as the first cell: $\Delta x_{-1} =
\Delta x_{0} = \Delta x_1$.

$\bullet$\ \textsc{Step 1}: computation of ${f}_{1/2,k}^{n,+}$ for
outgoing velocities $v_k <0$.\\
This value is defined by a linear extrapolation using the values
${f}_{1,k}^{n}$ and ${f}_{2,k}^{n}$: 
\begin{equation}\label{eq-extrap12} 
  {f}_{1/2,k}^{n,+} = \frac32{f}_{1,k}^{n} - \demi{f}_{2,k}^{n}.
\end{equation}

\smallskip

$\bullet$\ \textsc{Step 2}: computation of ${f}_{1/2,k}^{n,-}$  for incoming velocities $v_k >0$.\\
We use the boundary condition~(\ref{eq-BCdiffk_1d})--(\ref{eq-sigmadiffk_1d}) to get 
\begin{equation}
{f}_{1/2,k}^{n,-} = \sigma_w ({f}_{1/2}^{n,+}) \mathcal{M}^n_k.\label{diff2}
\end{equation}

\smallskip

$\bullet$\ \textsc{Step 3}: definition of the values ${f}_{0,k}^{n}$ and ${f}_{-1,k}^{n}$ in the ghost cells.\\

For outgoing velocities, only the value of $f^n_{0,k}$ is required. We
use the same linear extrapolation as used to compute extrapolation
${f}_{1/2,k}^{n,+}$ (see Step 1) to get
\begin{equation*}
    {f}_{0,k}^{n} = 2{f}_{1,k}^{n} - {f}_{2,k}^{n}.
\end{equation*}

For incoming velocities, we need the values of ${f}_{0,k}^{n}$ and
${f}_{-1,k}^{n}$. Then we use again a linear extrapolation, but based
on the values of the incoming value of the wall interface distribution
${f}_{1/2,k}^{n,-}$ and on $f^n_{1,k}$. This gives
\begin{align*}
  f^n_{0,k} & = 2 {f}_{1/2,k}^{n,-} - f^n_{1,k}, \\
  f^n_{-1,k} & = 4 {f}_{1/2,k}^{n,-} - 3 f^n_{1,k}. 
\end{align*}

Now it is not difficult to prove that these definitions satisfy the
previous constraints, provided that $\alpha_i\geq\demi$. First,~(\ref{eq-fdemicl}) is imposed at step 2,
and hence it is satisfied, which gives mass conservation. Now, it remains to prove that
$f^{n,\pm}_{\demi,k}$ defined by~(\ref{eq-extrap12}) and~(\ref{diff2}) also satisfy~(\ref{eq-fdemiplus})
and~(\ref{eq-fdemimoins}). This is due to the extrapolation procedures
that make the points $(x_i,f^n_{i,k})$ are on the same straight line,
and hence the limiters can be computed. Indeed, for instance, we have
for incoming velocities $v_k>0$: 
\begin{equation*}
\begin{split}
\delta_{0,k}^{n,lim} & = {\rm minmod}\left(
\frac{f^n_{1,k}-f^n_{-1,k}}{2\Delta x_0},
2 \alpha_0\frac{f^n_{0,k}-f^n_{-1,k}}{\Delta x_0}, 
2 \alpha_0 \frac{f^n_{1,k}-f^n_{0,k}}{\Delta x_0}
\right) \\
& = {\rm minmod}\left(
2\frac{f^n_{1,k} - f^{n,-}_{1/2,k}}{\Delta x_0},
\modif{(2 \alpha_0)2\frac{f^n_{1,k} - f^{n,-}_{1/2,k}}{\Delta x_0}, 
(2 \alpha_0)2 \frac{f^n_{1,k} - f^{n,-}_{1/2,k}}{\Delta x_0}}
\right) \\
& = 2\frac{f^n_{1,k} - f^{n,-}_{1/2,k}}{\Delta x_0}.
\end{split}
\end{equation*}

Consequently, the right-hand side of~(\ref{eq-fdemimoins}) is 
\begin{equation*}
  \begin{split}
  {{f}_{0,k}^n}+\frac{\Delta x_{0}}{2} {{\delta}_{0,k}^{n,lim}} & =
 2 {f}_{1/2,k}^{n,-} - f^n_{1,k} +  f^n_{1,k} - f^{n,-}_{1/2,k} \\
& = {f}_{1/2,k}^{n,-}
  \end{split}
\end{equation*}
and~(\ref{eq-fdemimoins}) is satisfied. For outgoing velocities, we
have
\begin{equation*}
\begin{split}
{\delta}_{1,k}^{n,lim} & = {\rm minmod}\left(
\frac{f^n_{2,k}-f^n_{0,k}}{2\Delta x_1},
2 \alpha_1\frac{f^n_{1,k}-f^n_{0,k}}{\Delta x_1}, 
2 \alpha_1 \frac{f^n_{2,k}-f^n_{1,k}}{\Delta x_1}
\right) \\
& = {\rm minmod}\left(
\frac{f^n_{2,k}-f^n_{1,k}}{\Delta x_1},
2 \alpha_1\frac{f^n_{2,k}-f^n_{1,k}}{\Delta x_1}, 
2 \alpha_1 \frac{f^n_{2,k}-f^n_{1,k}}{\Delta x_1}
\right) \\
& = \frac{f^n_{2,k}-f^n_{1,k}}{\Delta x_1}.
\end{split}
\end{equation*}
Consequently, the right-hand side of~(\ref{eq-fdemiplus}) is
\begin{equation*}
  \begin{split}
  {{f}_{1,k}^n}-\frac{\Delta x_{1}}{2} {{\delta}_{1,k}^{n,lim}}
&  = \frac32 {f}_{1,k}^n - \demi f^n_{2,k} \\
& = {f}_{1/2,k}^{n,+}
  \end{split}
\end{equation*}
from~(\ref{eq-extrap12}), and hence~(\ref{eq-fdemiplus}) is satisfied.

In summary, we have proved that this finite volume scheme
is based on slope limiters whose ghost cell values are consistent with
the slope reconstruction up to the solid boundary: this implies the
second order accuracy up to the boundary. Moreover, we have proved
that this scheme also preserves the mass conservation even with solid boundaries.
\textcolor{black}{Finally, note that even if the extrapolation procedure we use
  at the solid wall might induce non positive values of the
  distribution functions, we did not observe this problem in our
  numerical tests. This is probably due to the fact that we use small
  cells in the Knudsen layer, which make the gradients, and hence the
  slopes, small enough too avoid the creation of negative values.}

\begin{figure}[p]
\centering
\includegraphics[scale=0.5]{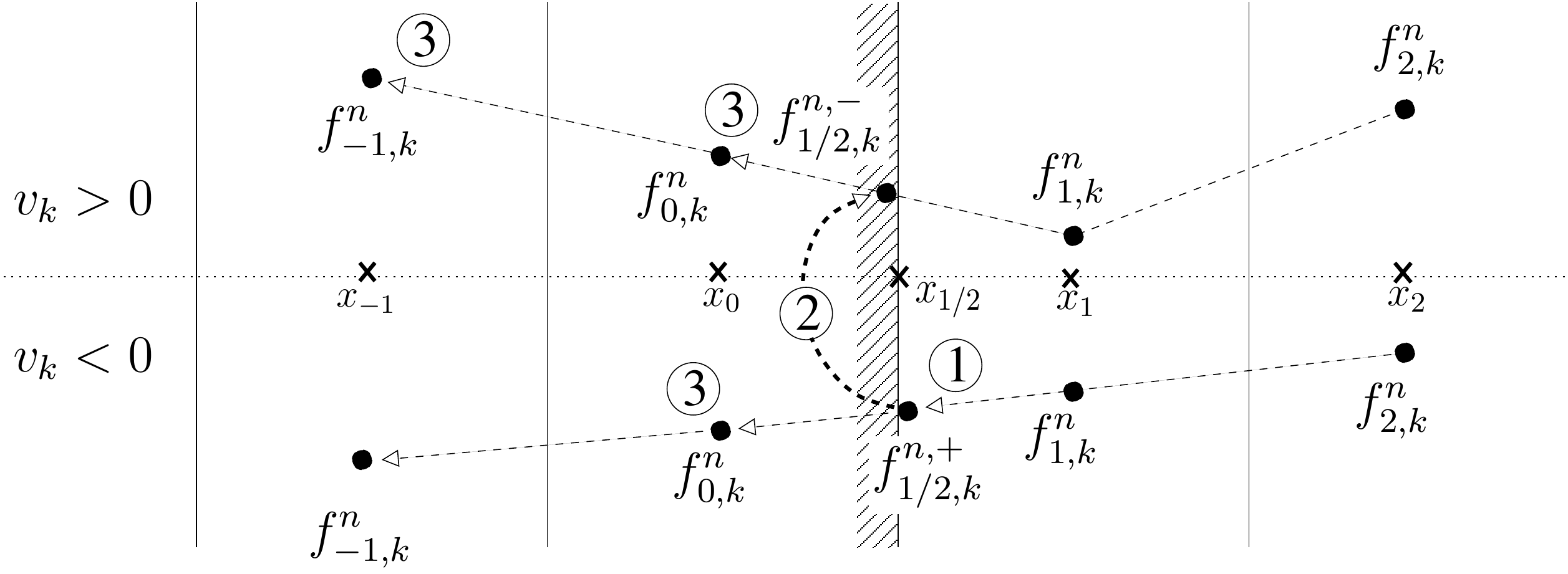}
\caption{\textcolor{black}{Finite volume scheme with the Yee flux limiter: definition of the ghost cell values
  $f^n_{0,k}$ and $f^n_{-1,k}$ by using linear extrapolations and the
  boundary condition.}}
\label{fig:O2flux}
\end{figure}

  \begin{remark}
\textcolor{black}{The fact that scheme with a flux limiter presented in
  section~\ref{schemyee} cannot be both
second order and conservative at the solid wall can be seen as
follows. 
}

\textcolor{black}{The previous analysis for the scheme with a slope
  limiter relies on the fact that the minmod
function~(\ref{limgeneral}) reduces to a single
slope. Consequently, the scheme is linear at the solid wall, and
hence is compatible with boundary condition~(\ref{eq-fdemicl})
(which is linear too). This implies the conservation property. This
reduction is due to the fact that the ghost cell values $f_{0,k}$ and
$f_{-1,k}$ are defined by extrapolations that make all the points used
in~(\ref{limgeneral}) {\it aligned} (and hence with three equal
slopes).
}

\textcolor{black}{For the scheme with a flux limiter, this extrapolation cannot make all
the points aligned. Indeed, the limiter~(\ref{eq-yee}) uses the
same stencil for negative and positive velocities (it is a
``symmetric'' limiter (see~\cite{yee}). This implies that for
positive velocities, $f^n_{-1,k}$ and $f^n_{0,k}$ are defined through
extrapolation of $f^{-,n}_{\demi,k}$ and $f^n_{1,k}$, but they cannot
be aligned with $f^n_{2,k}$, in general (see
figure~\ref{fig:O2flux}). Our previous analysis cannot be applied, and
there is no reason for the conservation property holds true.
}

\textcolor{black}{However, if instead a non symmetric (or ``upwind'') flux limiter is
used (see~\cite{Leveque_book}), it is still possible to prove the
conservation property. Such is scheme is not used in this paper.
}
  \end{remark}

\subsection{Discontinuous Galerkin scheme}
\label{subsec:DG}

In finite volume methods, the distribution is approximated by a
piecewise constant function (on each cell, the discrete distribution
is equal to its cell average). Discontinuous Galerkin schemes can be
viewed as an extension of the finite volume method in which the
distribution is now approximated by piecewise polynomial
functions. Here, we will consider a piecewise linear approximation.

\subsubsection{Weak form and piecewise linear approximation}
\label{subsec:weak}

On each cell $\Omega_i = [\ximdemi,\xipdemi]$, we use
two basis affine functions $\varphi_{i,1}$ and $\varphi_{i,2}$ such
that $\varphi_{i,1}(\ximdemi) = 1$, $\varphi_{i,1}(\xipdemi) = 0$, and
$\varphi_{i,2}(\ximdemi) = 0$, $\varphi_{i,2}(\xipdemi) = 1$. We
project equation~(\ref{BGKdiscret}) in this basis and integrate by
parts to get for every $i$ and $p=1,2$:
 \begin{equation} \label{wf}
\begin{split}
&  \int_{\Omega_i} \frac{\partial f_k}{\partial t}(t,x) \varphi_{i,p}(x) \, dx
 - \int_{\Omega_i} v_k f_k(t,x) \dx \varphi_{i,p}(x) \, dx \\
&  + \left(
v_kf_k(t,x_{i+\demi}) \varphi_{i,p}(x_{i+\demi}) 
- v_kf_k(t,x_{i-\demi}) \varphi_{i,p}(x_{i-\demi})
\right)
=  \int_{\Omega_i} \frac
{1}{\tau_i}\left(\mathcal{M}_{k}[\rho,\um,T]  - f_k\right)
\varphi_{i,p}(x) \, dx, 
\end{split}
\end{equation}

\begin{figure}[p]
\centering
\includegraphics[scale=0.5]{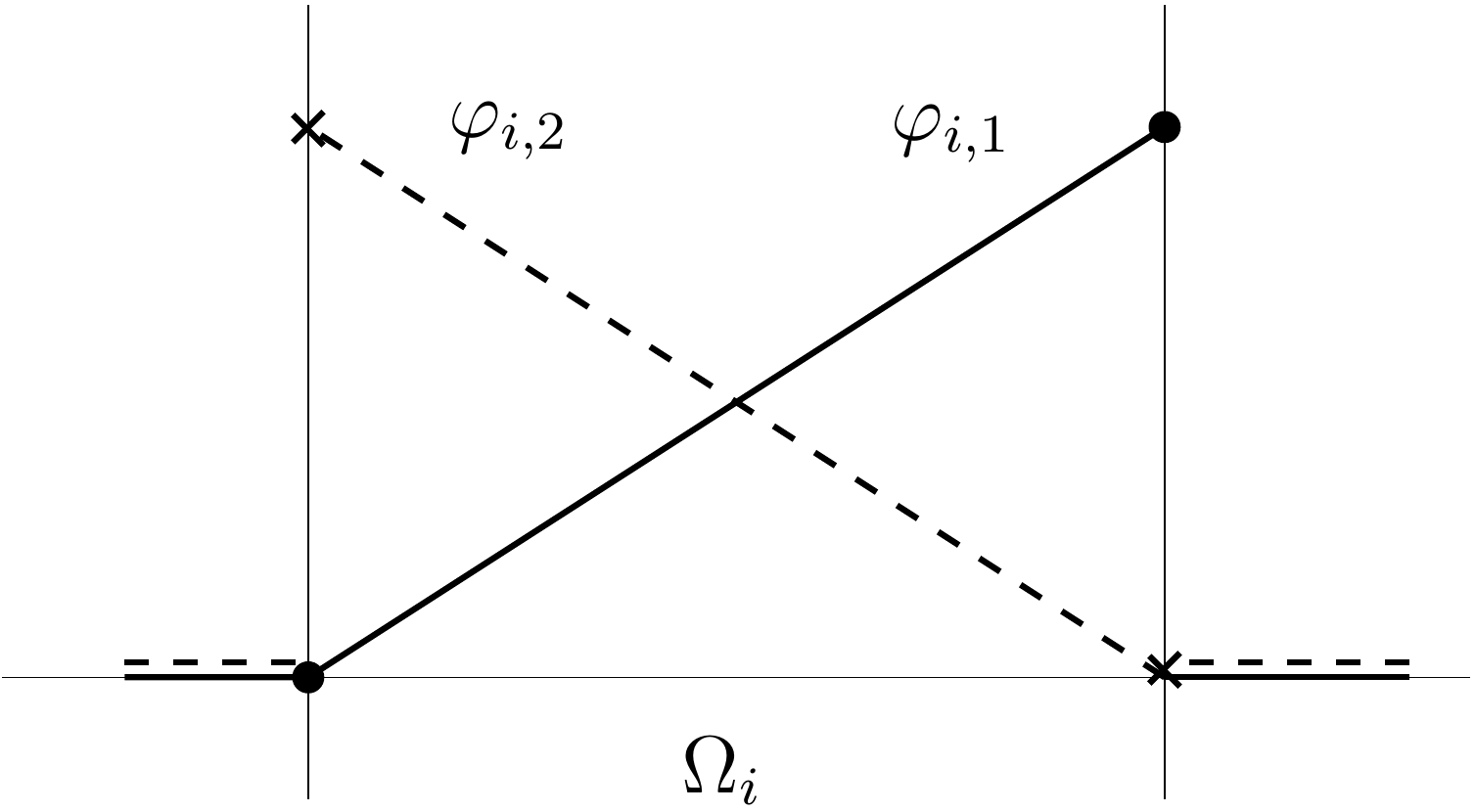} \hfill
\includegraphics[scale=0.5]{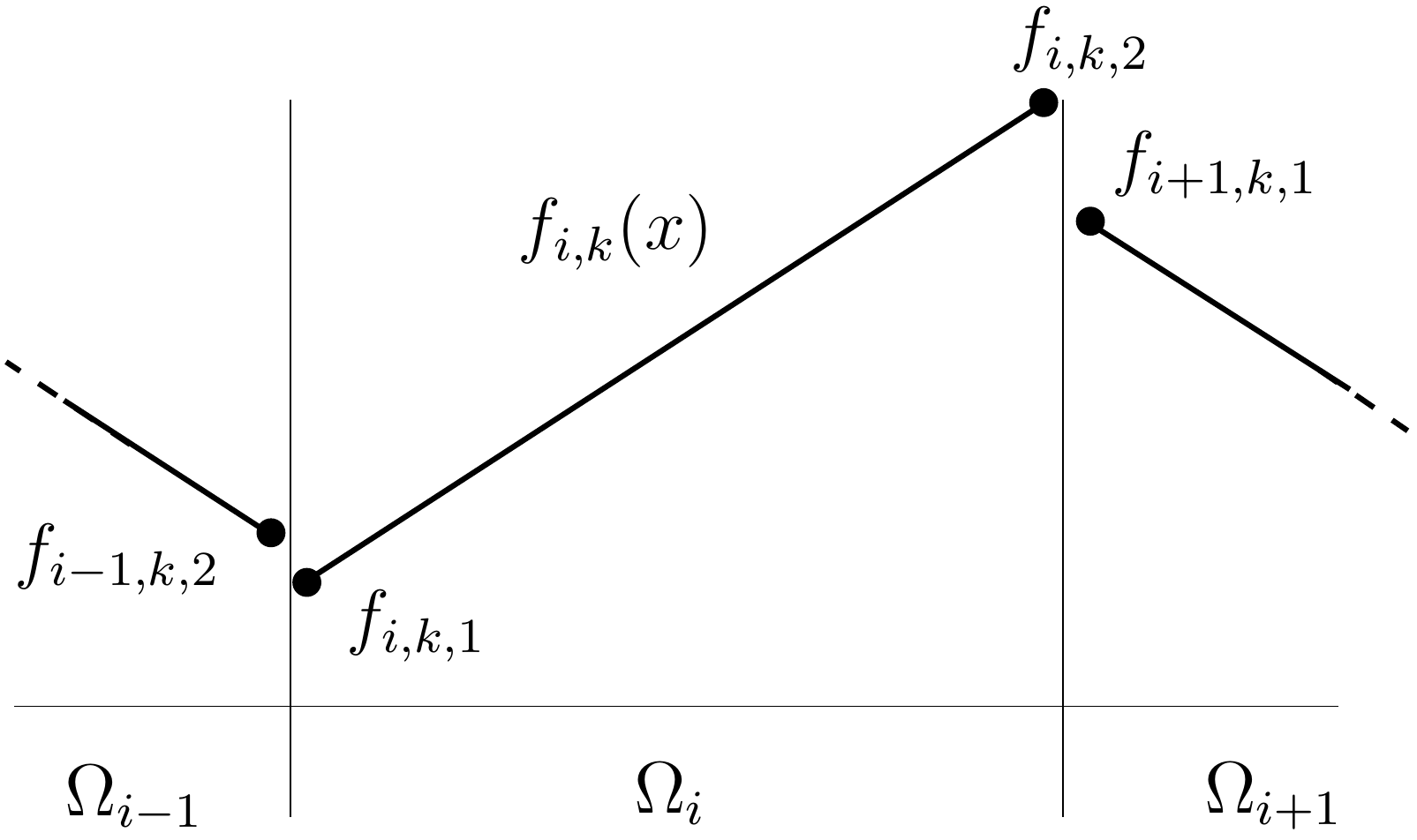} \\
\includegraphics[scale=0.5]{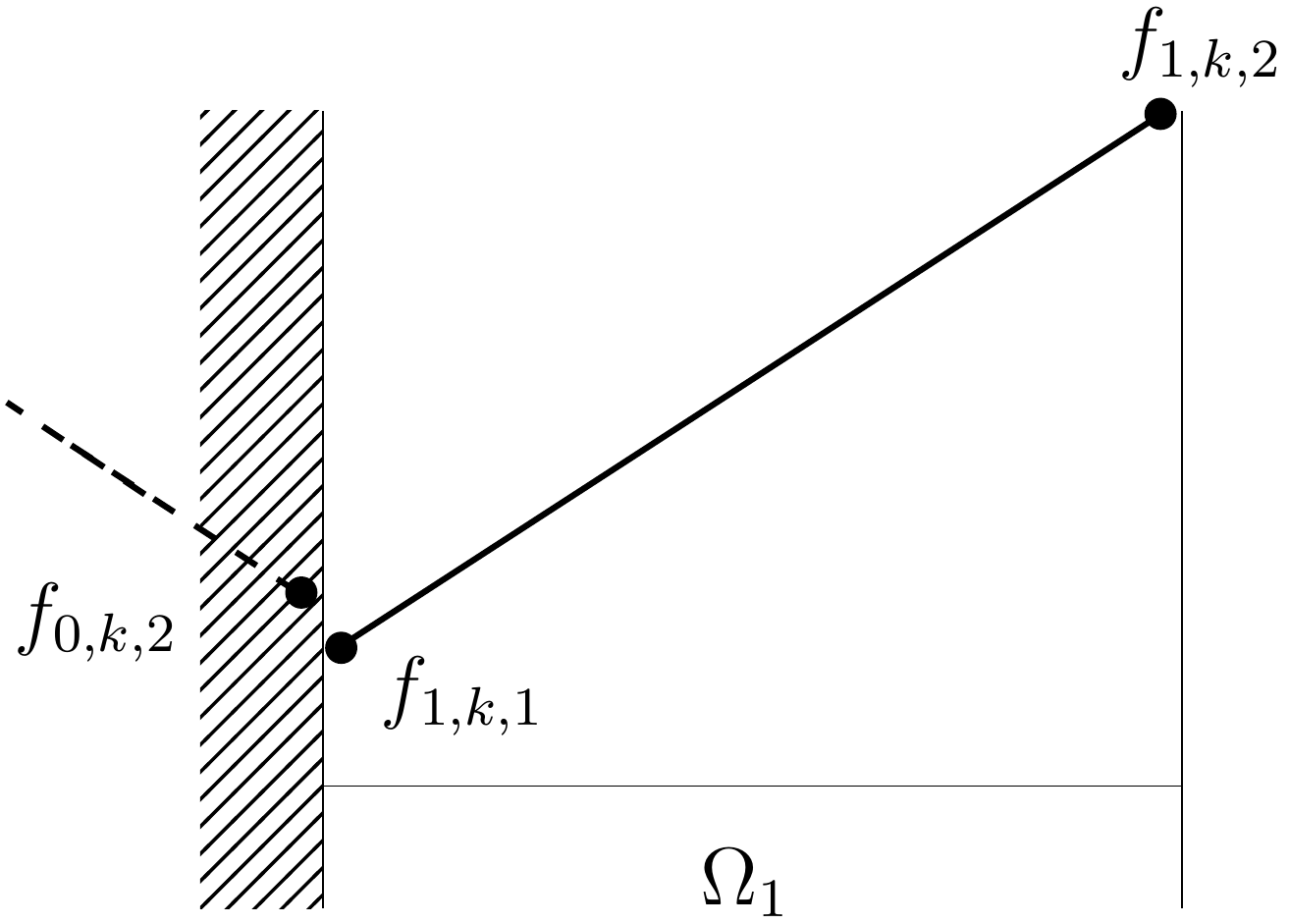}
\caption{Piecewise linear basis functions (left), piecewise
  linear approximation of $f_k$ (right), solution at the wall (bottom).}
\label{fig:DGapprox}
\end{figure}

Now, we assume that $f_k(t,x)$ can be approximated on each cell by a
piecewise linear function $\hat{f}_k(t,x)$ defined by 
$\hat{f}_k(t,x)|_{\Omega_i} = f_{i,k}(t,x) = 
f_{i,k,1}(t)\varphi_{i,1}(x) + f_{i,k,2}(t)\varphi_{i,2}(x) 
$ (see figure~\ref{fig:DGapprox}). Note that while the components
$f_{i,k,1}$ and $f_{i,k,2}$ are also the pointwise values of $f_{i,k}$
at the edges of $\Omega_i$, the whole function $\hat{f}_k$ itself is
not continuous across these edges: indeed, at each cell edge
$x_{i+\demi}$, $\hat{f}_k$ has two left and right values
$f_{i,k}(t,x_{i+\demi}) = f_{i,k,2}$ and $f_{i+1,k}(t,x_{x_{i+\demi}})
= f_{i+1,k,1}$ (see figure~\ref{fig:DGapprox}).

This approximation has now to be injected in the weak
form~(\ref{wf}). The first term of the left-hand side is easily computed
 \begin{equation}\label{eq-T1} 
\displaystyle  \sum_{q=1}^{2}{ \int_{\Omega_i}}  \frac {\partial
  f_{i,k,q}}{\partial t}\varphi_{i,q}(x)\varphi_{i,p}(x) \, dx 
= \sum_{q=1}^{2} m_{pq}^i \frac {\partial f_{i,k,q}}{\partial t},
 \end{equation}
where the $(m_{pq}^i)$ form the $2\times 2$ matrix 
\begin{equation}\label{eq-mpq} 
  (m_{pq}^i) = \bar{\bar{M_i}} = 
\frac{|\Omega_{i}|}{6}
\begin{pmatrix}
2 & 1 \\
1 & 2
\end{pmatrix}. 
\end{equation}

The second term of the left-hand side is
\begin{equation}\label{eq-T2} 
 -v_k\displaystyle{ \int_{\Omega_i}}f_{i,k}(x) \frac
{\partial \varphi_{i,p}(x)}{\partial x} \, dx 
  = \sum_{q=1}^{2} d_{pq}f_{i,k,q},
\end{equation}
where the $(d_{pq})$ form the $2\times 2$ matrix 
\begin{equation} \label{eq-dpq} 
(d_{pq}) = \bar{\bar{D}}_k  
= \frac{v_k}{2}
\begin{pmatrix}
 {\ \ 1} & {\ \ 1} \\
{-1} & {-1} \\
\end{pmatrix}.
\end{equation}

For the third term of the left-hand side of~(\ref{wf}), we muse take
into account that the piecewise linear approximation of $f_k$ is not
continuous across cell edges (see the remark above). Consequently, we
write this term as $v_k\hat{f}_k(t,x_{i+\demi}) \varphi_{i,p}(x_{i+\demi}) 
- v_k\hat{f}_k(t,x_{i-\demi}) \varphi_{i,p}(x_{i-\demi})$, where a
single value for $\hat{f}_k(t,x_{i+\demi})$ has to be defined. A good
choice (for accuracy, stability, and treatment of boundary conditions)
is to use the upwind value
\begin{equation*}
  \hat{f}_k(t,x_{i+\demi})  = 
  \begin{cases}
       f_{i,k,2}(t)      & \text{if $v_k>0$}, \\
 f_{i+1,k,1}(t) & \text{if $v_k<0$}.
  \end{cases}
\end{equation*}
Taking into account that the basis functions $\varphi_{i,p}$ take
values $0$ or $1$ at the cell edges of $\Omega_i$, we get
\begin{equation}  \label{eq-T3}
v_k\hat{f}_k(t,x_{i+\demi}) \varphi_{i,p}(x_{i+\demi}) 
- v_k\hat{f}_k(t,x_{i-\demi}) \varphi_{i,p}(x_{i-\demi}) 
=
\begin{cases}
   -v_k^+ f_{i-1,k,2}(t) - v_k^- f_{i,k,1}(t) & \text{for $p=1$}, \\
 v_k^+ f_{i,k,2}(t) + v_k^- f_{i+1,k,1}(t) & \text{for $p=2$}.
\end{cases}
\end{equation}

The right-hand side of~(\ref{wf}) is obtained by using the same
piecewise linear approximation for each term, which gives
\begin{equation}  \label{eq-T4}
\int_{C_i}\frac {1}{\tau}(\mathcal{M}_k[\hat{\textbf{U}}]  -
\hat{f}_{k})\varphi_{i,p}(x) \, dx 
=\sum_{q=1}^{2}m_{pq}^i\frac {1}{\tau_{i,q}}\left(\mathcal{M}_{i,k,q}-f_{i,k,q}\right),
\end{equation}
where $m_{pq}^i$ has been defined in~(\ref{eq-mpq}),
$\mathcal{M}_{i,k,q} = \mathcal{M}_k[\rho_{i,q},\um_{i,q},T_{i,q}]$, 
$\tau_{i,q} = \tau(\rho_{i,q},T_{i,q})$, and the macroscopic 
quantities $\rho_{i,q},\um_{i,q}, T_{i,q}$ are defined by applying~(\ref{eq-mtsk}) in each cell.

Finally, we collect the different terms~(\ref{eq-T1}--\ref{eq-T4}),
and we find that the approximation of~(\ref{wf}) can be written in the
following vectorial form 
 \begin{equation}\label{eq-vecformDG} 
\bar{\bar{M_i}} \  \frac {\partial {\bf f}_{i,k} }{\partial
  t}+\bar{\bar{D}}_k \ {\bf f}_{i,k} + 
\left( 
\bar{\bar{A}}_k \ {\bf f}_{i-1,k}
+ \bar{\bar{B}}_k \ {\bf f}_{i,k}
+ \bar{\bar{C}}_k \ {\bf f}_{i+1,k}
  \right)
= 
 \bar{\bar{M_i}} \ \frac{1}{\boldsymbol{\tau}_i} ({\boldsymbol {\cal {M}}}_{i,k}-{\bf f}_{i,k}), 
 \end{equation}
where we use the following notations for the
two-component vectors ${\bf f}_{i,k} = (f_{i,k,1},f_{i,k,2})^T$ and
${\boldsymbol {\cal {M}}}_{i,k} =
(\mathcal{M}_{i,k,1},\mathcal{M}_{i,k,2})^T$, the 2$\times$2 matrices $\bar{\bar{A}}_k =
\left(\begin{smallmatrix}
  0 & -v_k^+ \\ 0 & 0
 \end{smallmatrix}\right)
$, 
$\bar{\bar{B}}_k =
\left(\begin{smallmatrix}
  -v_k^- & 0 \\ 0 & v_k^+
 \end{smallmatrix}\right)
$, 
$\bar{\bar{C}}_k =
\left(\begin{smallmatrix}
  0 & 0 \\ v_k^- & 0
 \end{smallmatrix}\right)
$, while the matrices $\bar{\bar{M_i}}$ and $\bar{\bar{D}}_k$ have been
defined in~(\ref{eq-mpq}) and~(\ref{eq-dpq}).
Note that the product of vectors $\frac{1}{\boldsymbol{\tau}_i} ({\boldsymbol {\cal
    {M}}}_{i,k}-{\bf f}_{i,k})$ must be understood component-wise.
 Multiplying~(\ref{eq-vecformDG}) by
the matrix ${\bar{\bar{M_i}}}^{-1} ={\frac{2}{|\Omega_i|}}
\left(\begin{smallmatrix}
    \phantom{-}2 & -1 \\
    -1 & \phantom{-}2
 \end{smallmatrix}\right)
$ 
(the inverse of
$\bar{\bar{M_i}}$), we obtain the semi-discrete scheme
 \begin{equation}
\begin{aligned}
\frac {\partial {\bf f}_{i,k} }{\partial t}=
&- \frac{1}{|\Omega_i|}
 \begin{pmatrix}
 {\ \ 3 v_k} & {\ \ 3 v_k} \\
 {-3 v_k} & {-3 v_k} \\
 \end{pmatrix} 
{\bf f}_{i,k}  \\
&- \frac{2}{|\Omega_i|} \left[
\begin{pmatrix}
\ \ 0\ \  & -2v_k^{+} \\
\ \ 0\ \ & v_k^{+}
\end{pmatrix} 
{\bf f}_{i-1,k} 
+
 \begin{pmatrix}
-2 v_k^{-} & -v_k^{+} \\
 v_k^{-}& 2 v_k^{+}
 \end{pmatrix} 
{\bf f}_{i,k} 
+
\begin{pmatrix}
-v_k^{-}& \ \ \ 0\ \   \\
2 v_k^{-}& \ \ \ 0\ \  
\end{pmatrix} 
{\bf f}_{i+1,k} 
\right] \\
&+
\frac{1}{\boldsymbol{\tau}_i}({\boldsymbol {\cal {M}}}_{i,k}-{\bf f}_{i,k}). \\
\label{gdsemidisc}
\end{aligned}
 \end{equation}

\subsubsection{Time discretization }
\label{subsec:time-disc}
As shown in \cite{shu}, Discontinuous Galerkin scheme are unstable
when used with a forward Euler method. Since we are interested in
steady flows only, we use the same linearized implicit scheme as for
the finite volume method (see~\cite{luc_jcp,Baranger}). 


\subsubsection{Boundary conditions}
\label{subsec:bound_conds}
The Discontinuous Galerkin scheme is very compact, as the value ${\bf
  f}_{i,k}^{n+1}$ only depends on  ${\bf f}_{i,k}^{n}$,
$ f_{i-1,k,2}^{n}$ and $f_{i+1,k,1}^{n}$. We just need one rank of
ghost cells on each boundary. Moreover, the value of $f_k^{n}$ in the
ghost cell has just to be defined on the side corresponding to the
interface between the wall and the gas, this means  $f_{0,k,2}^{n}$
for the left boundary (and $f_{i_{max}+1,k,1}^{n}$ at the right one).

For the left wall, we use a technique which is close to what is done
with the finite volume scheme (except that we do not need
extrapolation here): we define the distribution at the wall by 
\begin{equation*}
  f^n_{wall,k} =
  \begin{cases}
    f^n_{1,k,1} & \text{for $v_k<0$},\\
f^n_{0,k,2} & \text{for $v_k>0$},
  \end{cases}
\end{equation*}
see figure~\ref{fig:DGapprox}. The ghost cell value $f^n_{0,k,2}$ is
defined by applying the boundary condition~(\ref{eq-BCdiffk_1d})--(\ref{eq-sigmadiffk_1d}) to
$f^n_{wall,k}$ to get for $v_k>0$
\begin{equation*}
  f_{0,k,2} =  \sigma_1^{n} \ \mathcal{M}_k[1,u_w,T_w],
\end{equation*}
where $
\sigma_1^{n}  = -\left({\sum_{k}}v_k^{-}
    {f}_{1,k,1}^{n}\omega_k\right) / \left(\sum_{k}v_k^{+}
    \mathcal{M}^n_k[1,u_w,T_w]\omega_k\right)$.


\section{Numerical results}
\label{sec:Num}
 \subsection{Test-case and reduced model in 1D: Couette flow}
      \label{subsec:couette}

We consider a plane Couette flow of argon: the gas lies between two
flat walls, both walls are at the same temperature $T_w$, the left
wall is at rest, while the right wall moves upward with the velocity
${\bf {u}_{w}} = (0,u_w,0)^T$. The resulting flow is one
dimensional. With the standard reduced distribution
technique~\cite{chuchu}, the 3D BGK equation~(\ref{eq-kin}) reduces to
the following system of 3 kinetic equations
\begin{equation*}
  \dt
  \begin{pmatrix}
    F \\ G \\ H
  \end{pmatrix}
+ v_x \dx
\begin{pmatrix}
  F\\ G \\ H
\end{pmatrix} = 
\frac{1}{\tau}
\begin{pmatrix}
  \mathcal{M}[\rho,\bf{u},T] - F \\
  \mathcal{N}[\rho,\bf{u},T] - G \\
  \mathcal{P}[\rho,\bf{u},T] - H \\
\end{pmatrix},
\end{equation*}
where the reduced distributions are
\begin{equation*}
  \begin{pmatrix}
        F \\ G \\ H
  \end{pmatrix} 
= \int_{\R^2}
\begin{pmatrix}
  1 \\\demi (v_x^2 +v_y^2) \\ v_z
\end{pmatrix}
f \, dv_x dv_y,
\end{equation*}
the macroscopic quantities mass density, velocity, temperature, and
normal heat flux, are
\begin{equation*}
  \begin{pmatrix}
    \rho \\ \rho u_x \\ \rho u_y \\ \frac{3}{2}\rho R T 
  \end{pmatrix}
 = \int_{\R^2}
 \begin{pmatrix}
   F \\ v_x F \\ H \\ \demi (v_x-u_x)^2 F + G
 \end{pmatrix}
\, dv_x,
\end{equation*}
and
\begin{equation*}
  q_x = \int_{\R^2} \demi(v_x-u_x)^3 F + (v_x - u_x) (G-u_yH) \, dv_x,
\end{equation*}
and the reduced Maxwellians are
\begin{gather*}
  \mathcal{M}[\rho,\um,T]:= \frac{\rho}{\sqrt{2\pi RT}}\exp \left(-\frac{|v_x-u_x|^2}{2RT}\right), \\
  \mathcal{N}[\rho,\um,T]:=  \left(RT + \frac {u_y^2}{2}\right)
\mathcal{M}[\rho,\um,T], \quad \text{ and } \quad
\mathcal{P}[\rho,\um,T]:=  u_y \ \mathcal{M}[\rho,\um,T].
\end{gather*}
The diffuse reflection boundary condition
\eqref{eq-BCdiff}--\eqref{eq-sigmadiff} reduces to
\begin{equation*}
  \begin{pmatrix}
    F \\ G \\ H
  \end{pmatrix}
|_{\text{wall},\ \textbf{v} \cdot \textbf{n} > 0} = 
\sigma_w(F)
\begin{pmatrix}
  \mathcal{M}[1,\um_w,T_w] \\ \mathcal{N}[1,\um_w,T_w] \\ \mathcal{P}[1,\um_w,T_w]
\end{pmatrix},
\end{equation*}
where $\sigma_w(F) = -\left(\int_{v_x n(x) < 0}v_xF(t,\x,v_x) 
  \, dv_x \right)\left(\int_{v_x n(x)>0}
    v_x\mathcal{M}[1,\boldsymbol{u}_w,T](v_x)  \, dv_x \right)^{-1},$
with $n(x)=1$ at the left wall and $n(x) = -1$ at the right wall.

All the schemes described in the previous section readily apply to
this system.

 \subsection{Comparison of second order Finite Volume schemes on the
   1D Couette flow}
   \label{subsec:compVF1D}
   
   
The flow parameters and gas properties are given in Table
$\ref{couette}$. The Knudsen number is based on the distance between
the plates, and the mean free path is defined by (see Bird \cite{bird}, with the tables for the values of $\mu_0$, $T_0$ and $\omega$) 
\begin{equation*}
\lambda = \displaystyle{\left[\sqrt{2}\pi d^2 \frac{\rho}{m}\left(\frac{T_0}{T}\right)^{\omega-\frac{1}{2}}\right]^{-1}}
\ \ \textrm{with}\ \ 
d = \displaystyle{\left[\frac{5(\alpha+1)(\alpha+2)\sqrt{{\pi}^{-1}m k_b T}}{4 \alpha (5-2\omega)(7-2\omega) \mu_0}\right]^{\frac{1}{2}}}
\end{equation*}\\
where $k_b$ is the Boltzmann constant and $m$ the molecular mass of
the gas.

For this problem, we use a uniform mesh of 100
cells. The discrete velocity grid is a uniform grid of 40 points with bounds $\pm 953 \ \textrm{m.s}^{-1} $. 
We compare here four schemes:
\begin{itemize}
\item the second order finite volume scheme with the Yee flux
  limiter~(\ref{eq-yee}) and first order boundary conditions, as
  described in section~\ref{schemyee}, and denoted here by ``\yee'';
\item the second order finite volume scheme using linear
  reconstruction~(\ref{eq-rec}),  with second order boundary
  conditions~(\ref{eq-fdemiplus}) and~(\ref{eq-fdemimoins}), without
  limiter (i.e. with slope~(\ref{moindrecarre})), denoted by
  ``\slopenolim'';
\item the same scheme but with limited slope~(\ref{limgeneral}) with $\alpha_i =
  \displaystyle{1/2}$, denoted by ``\slope'';
\item the same scheme but with boundary conditions of first order
  only, denoted by ``\slopeCLzero''
  (that is to say with distributions at the wall defined
  by~(\ref{eq-fdemiplus}) and~(\ref{eq-fdemimoins}) with zero limited slopes).
\end{itemize}

\begin{figure}[p]
\centering
\includegraphics[scale=0.35]{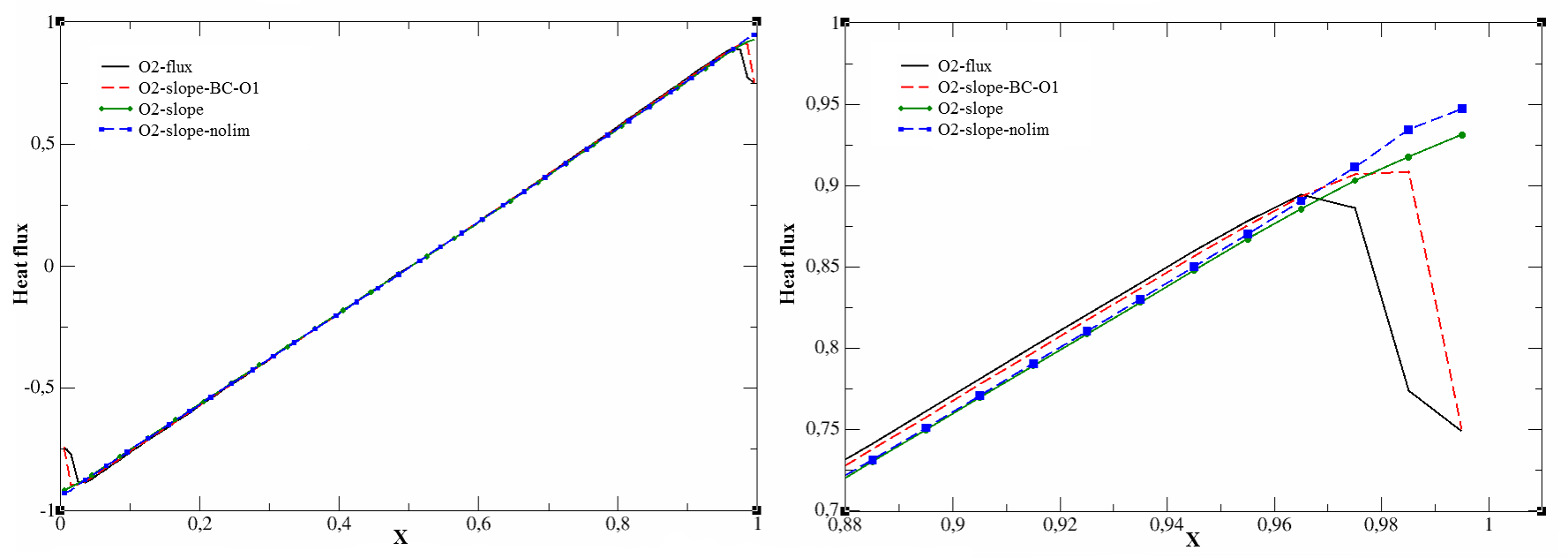}
\includegraphics[scale=0.35]{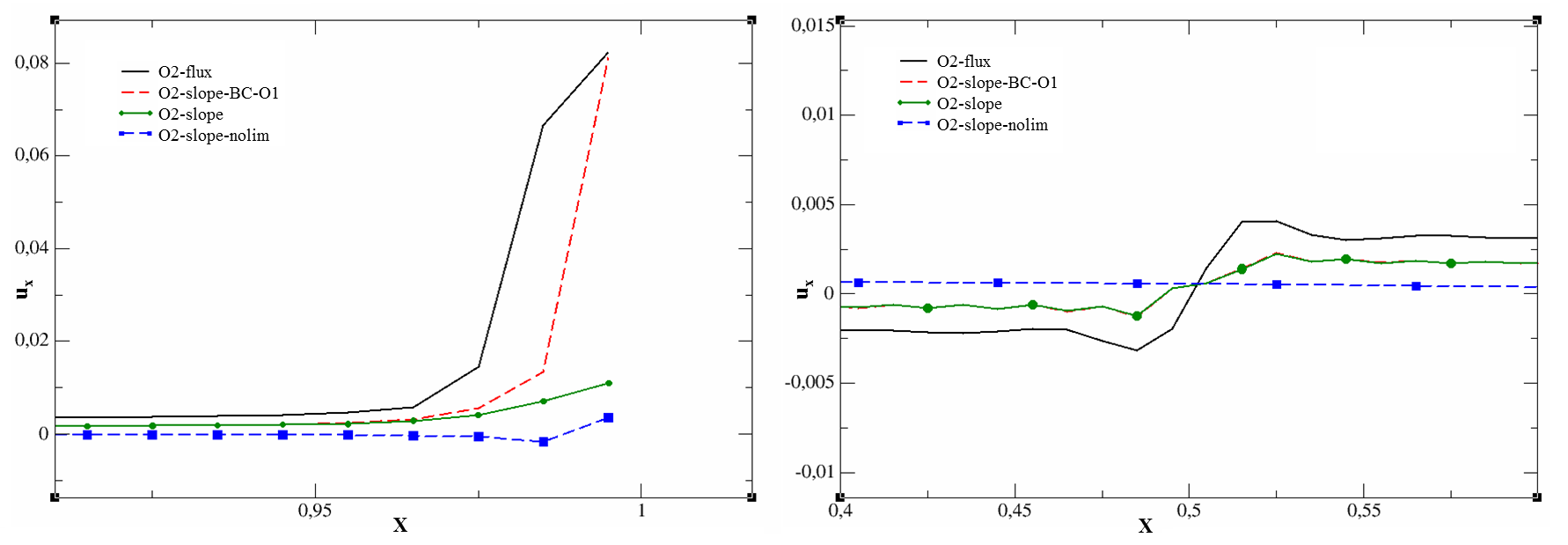}
\caption{1D Couette flow: {heat flux and horizontal velocity obtained with the
    different second order schemes: (top left) heat flux, (top right)
    heat flux, zoom at the moving wall, (bottom left) horizontal
    velocity, zoom at the moving wall, (bottom right) horizontal
    velocity, zoom in the middle of the domain.}}
\label{comp}
\end{figure}

We show on Figure $\ref{comp}$ the heat flux and the horizontal
velocity obtained with these four schemes. Note that the horizontal
velocity is a good indicator of accuracy, since its exact value is
zero in the whole domain.
We can see that the schemes with first order
boundary conditions (\yee  \ and \slopeCLzero) have a very bad behavior in the
near-wall zone, where we observe a strong variation of the velocity
and the heat flux, while the two other second order schemes give a
correct solution. This result confirms the analysis of section
\ref{schemyee}: the flux-limiter scheme degenerates into a first order
scheme in the near-wall zone, and the boundary conditions have to be
discretized in a way which is consistent with the second
order scheme.

Note that the slope limitation does not seem to be useful here: the scheme
without limitation is the most accurate, as this can be seen from
the horizontal velocity profiles (its exact value is zero). However,
we will see in next section that this limitation is
necessary for multidimensional cases.

   \subsection{Validation of the Finite Volume scheme in 2D}
\label{subsec:2D}
   The kinetic code of the CEA is dedicated to the simulation of
   hypersonic flows based on the BGK model, to make 2D plane,
   axisymmetric and 3D simulations of rarefied flows for monoatomic
   and diatomic gases. We refer to \cite{Baranger} for more details on
   this code. It uses two second order finite volume schemes on
   curvilinear meshes. The first one is the version with flux limiter method,
   as described in section~\ref{schemyee}, with first order boundary
   conditions. The other scheme is the one with the linear reconstruction and
   second order boundary conditions, as described in
   section~\ref{reconslin}. Since the space mesh is not uniform,
   formulas in section \ref{reconslin} (in particular
   $\eqref{moindrecarre}$) have to be modified accordingly.

\begin{figure}[p]
\centering
\includegraphics[scale=0.5]{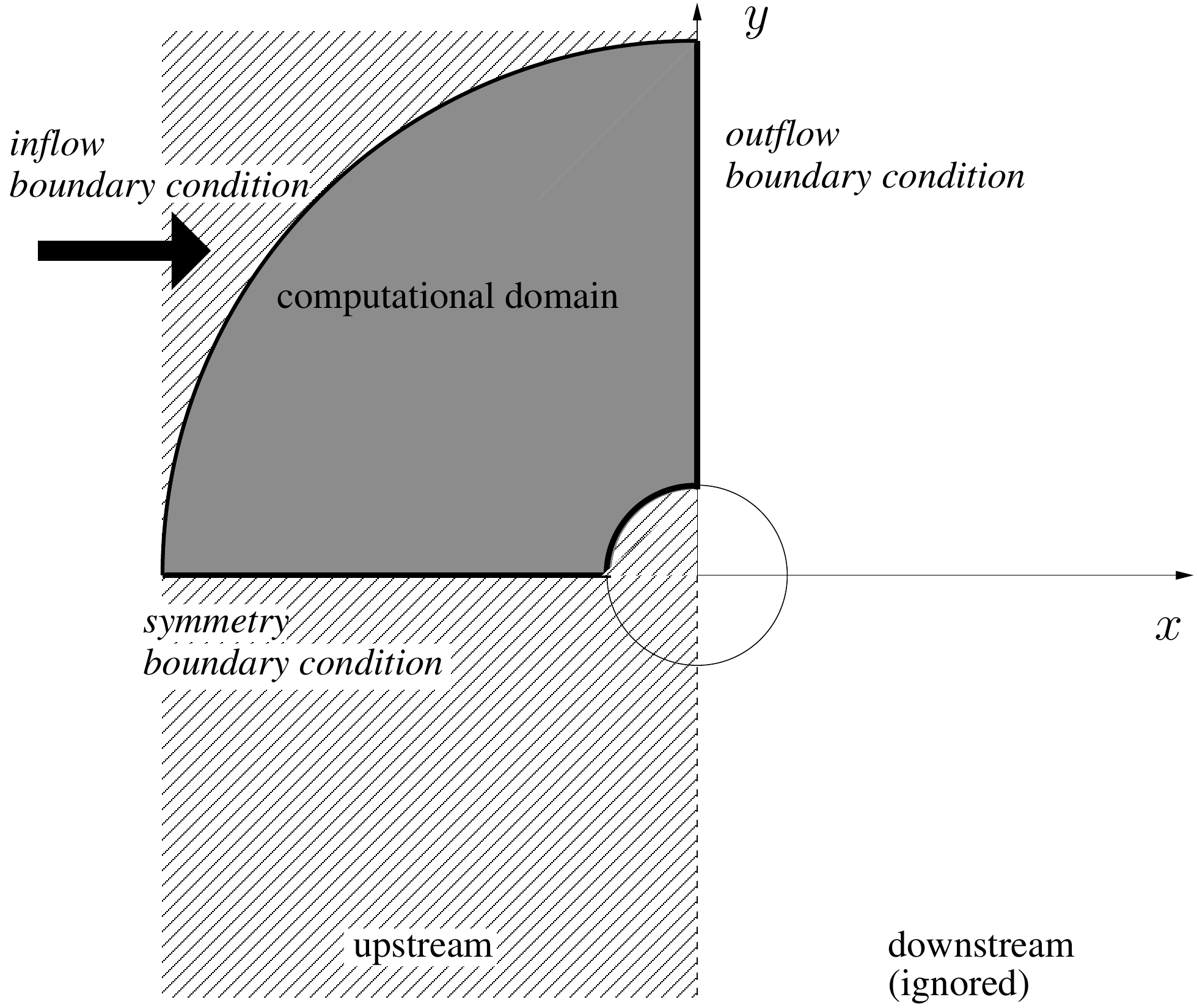}
\caption{\modif{2D supersonic flow around a cylinder: geometry of the test case.}}
\label{fig:geom_cyl}
\end{figure}

\begin{figure}[p]
\centering
\includegraphics[scale=0.5]{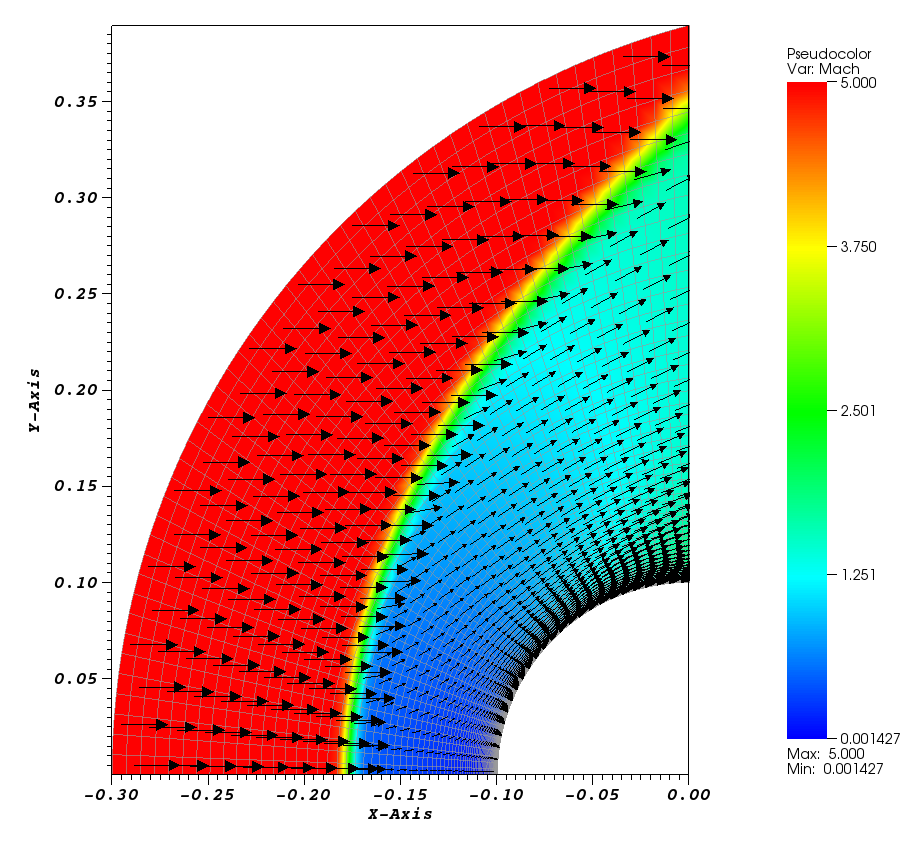}
\caption{\modif{2D supersonic flow: velocity and Mach number fields
    (finite volume scheme with slope limiters). The mesh used for the
    simulation is shown in grey lines.}}
\label{fig:cylindre_champs}
\end{figure}

   The test case is a steady flow of argon over a cylinder of radius
   0.1 m at Mach 5, \modif{see figure~\ref{fig:geom_cyl}}. The density
   and pressure are that of the air at an altitude of 60 km. Namely,
   we have $\rho = 3.059\ 10^{-4}\textrm{kg/m}^3$,
   $u = 1750\ \textrm{m/s}$ and $T = 352.6\ \textrm{K}$. The Knudsen
   number, based on the radius of the cylinder, is
   $\textrm{Kn}= 2.2\ 10^{-3}$, so the flow is rather dense, and a
   comparison with a Navier-Stokes simulation (with no-slip and no
   temperature jump) is relevant. With such a small Knudsen number, we
   can expect a difference between BGK and Navier-Stokes results lower
   than 1\%. \modif{As an illustration, the 2D velocity and Mach number
   fields are shown in figure~\ref{fig:cylindre_champs}, with the mesh
 used for the simulation, made with the finite volume scheme with
 slope limiters.}

  We compute the normal component of the heat flux along
   the boundary of the cylinder and
   compare the following finite volume schemes:
   \begin {itemize}
   \item first order scheme (\first);
   \item second order scheme with Yee flux limiters 
with first order boundary conditions (\yee)
   \item second order scheme with slope limiters and second order
     boundary conditions (\slope);
 \item Navier-stokes solution (\NS),  obtained with a finite volume code
   of the CEA. 
   \end{itemize}

\begin{figure}[p]
\centering
\includegraphics[scale=0.5]{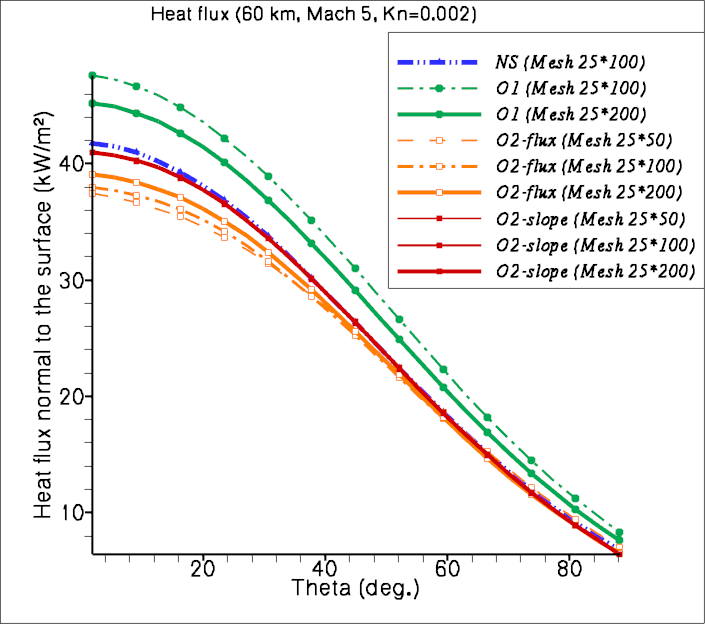}
\caption{2D supersonic flow: {comparison of the heat flux between BGK solved with
    finite volume schemes (\first, \yee, \slope) for different meshes, and Navier-Stokes
    (\NS).}\modif{ Profile of the normal heat flux along the solid
    wall.}}
\label{Comp2D}
\end{figure}

   In figure~\ref{Comp2D}, we show a comparison of the normal
   component of the heat flux along the solid boundary obtained with
   these schemes on different meshes (from $25\times 50$ to $25\times
   200$ cells), with an increasing number of cells in the direction
   orthogonal to the solid boundary. Like in the 1D Couette flow, we
   observe that the second order scheme with flux limiter is not
   sufficient to compute accurately the heat flux, while the second
   order scheme with linear reconstruction, slope limiters, and second
   order boundary conditions, is much more accurate.  Note that for
   this test case, this last scheme requires a real limitation: the
   coefficient $\alpha_i$ of the MC limiter~(\ref{limgeneral}) has
   been taken equal to $0.75$. \textcolor{black}{Indeed, the code
     produces too strong oscillations for
   larger values: this induces negative density and temperature that
   make the code stop.}

\subsection{Comparison between the Discontinuous Galerkin and Finite Volume schemes}
\label{subsec:compDGVF}


We compare the Discontinuous Galerkin scheme presented in
section~\ref{subsec:DG} to the finite volume schemes presented in
section~\ref{subsec:VF} and already analyzed in the numerical
comparisons above.

\begin{figure}[p]
\centering
\includegraphics[width=\textwidth]{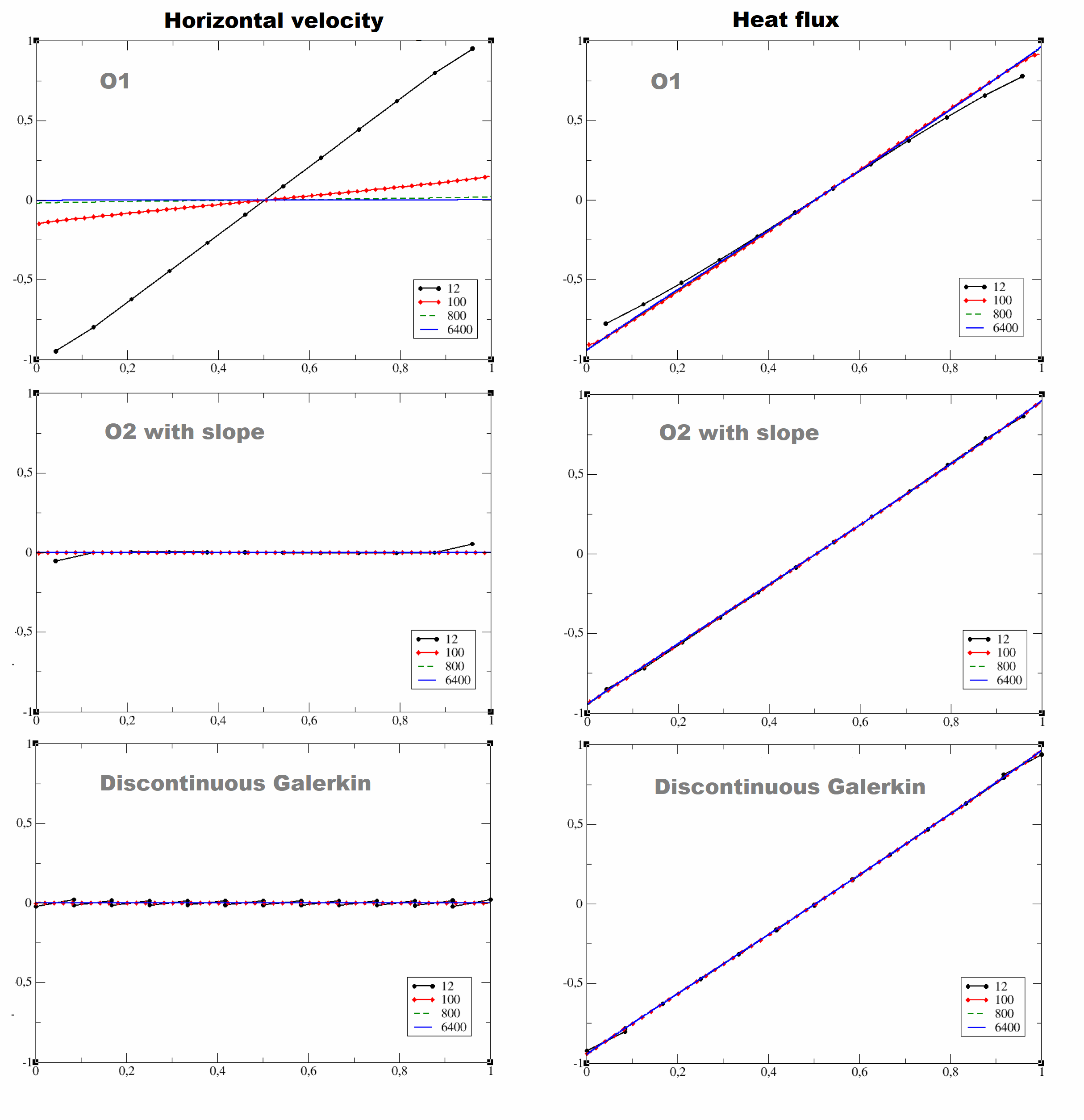}
\caption{1D Couette flow: results obtained with the discontinuous Galerkin
    scheme and first and second order finite volume schemes for the
    horizontal velocity and the heat flux, with 12 to 6400 grid
    points.}
\label{comp1}
\end{figure}

The test case is the same Couette flow as in section~\ref{subsec:compVF1D}.
We use uniform meshes with 12, 100, 800 and 6400 cells. We present on
Figure $\ref{comp1}$ the horizontal velocity (which should be exactly
zero) and the heat flux obtained for each scheme on each mesh. We can
observe that the Discontinuous Galerkin scheme is very close to the
second-order finite volume one; the results are far more accurate than
with the first-order finite volume scheme (note that the second order
scheme with flux limiters is not shown here).

With those results, we can compute the convergence rate of each
scheme, given in Table~\ref{tabconv} for the temperature and the
heat flux\modif{, in the $L^2$ norm}. 
\begin{table}[h]
\begin{center}
\begin{tabular}{|c|c|c|}
\hline
      & $q_x$ & $T$ \\ \hline
\first   & 0.98 & 1.06\\ \hline
\slope  & 1.89& 1.96\\ \hline
\yee   & 1.44& 1.84\\ \hline
\DG  & 1.74& 1.64\\ \hline
\end{tabular}
\end{center}
\caption{{\color{black}1D Couette flow: convergence orders in $L^2$ norm for the first order upwind 
finite volume scheme (\first), second order finite volume scheme with
slope limiters (\slope), finite volume scheme with flux limiters
(\yee), discontinuous Galerkin scheme (\DG)}.}
\label{tabconv}
\end{table} 
\modif{The convergence rates are consistent with our analysis and with the
theory. The upwind finite volume scheme shows a first order
convergence rate, the finite volume scheme with slope limiters shows an
almost second order convergence rate, the discontinuous Galerkin
scheme with upwind fluxes shows an almost 3/2 order convergence rate
(see~\cite{JP_1986}), while, as expected, the finite volume scheme with
flux limiters has an order of convergence smaller than 2.}

\begin{figure}[p]
\centering
\includegraphics[scale=0.45]{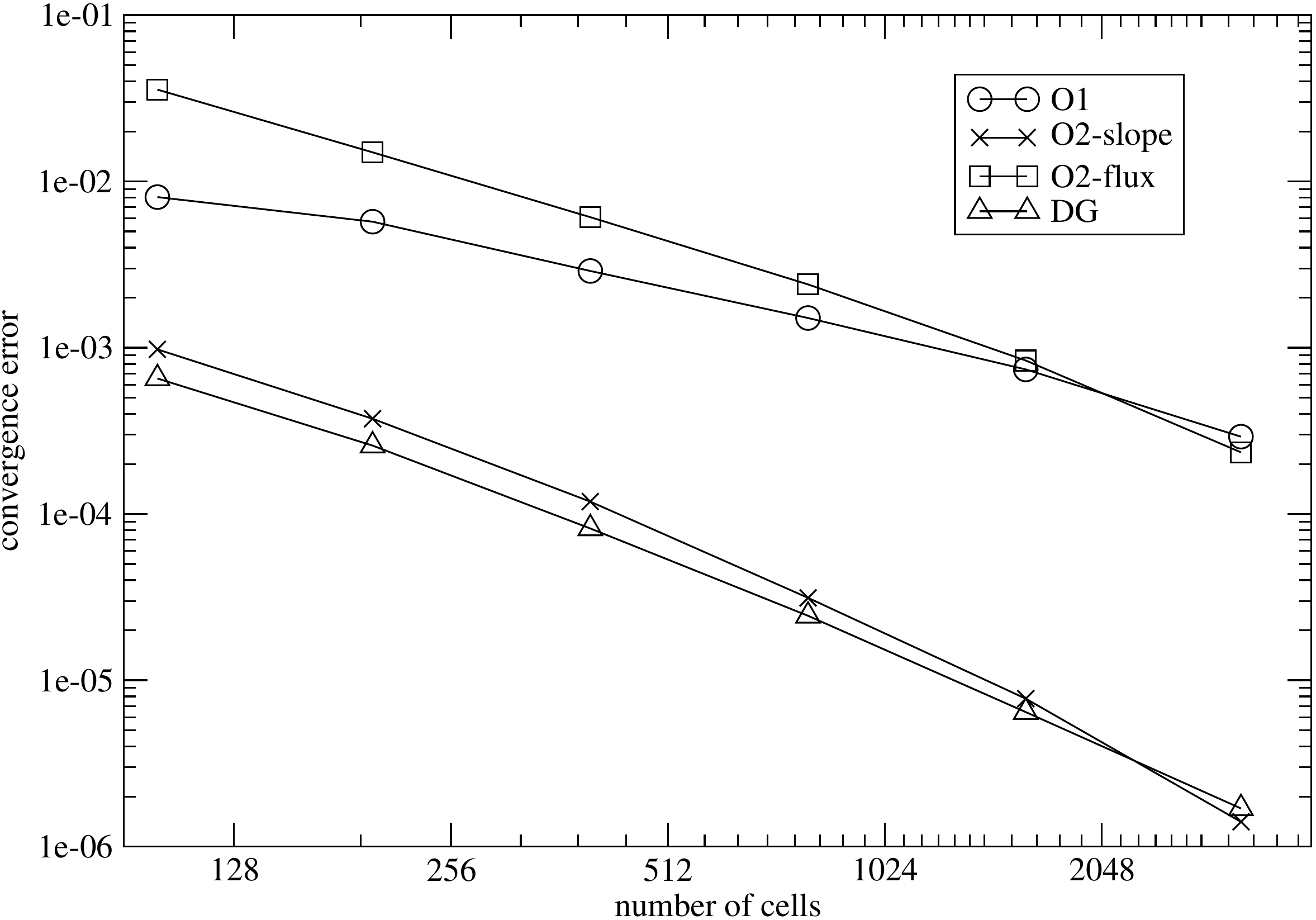}
\caption{1D Couette flow: 
\modif{mesh convergence ($L^2$ norm)} of the discontinuous Galerkin and
    finite volume schemes with uniform meshes, for the heat flux.
} 
\label{courbesconv}
\end{figure}

\modif{The corresponding error curves for the heat flux are given in
Figure~\ref{courbesconv}. It can be seen that the finite volume scheme with
flux limiters is not very accurate: even the first order scheme is
more accurate, while the numerical boundary flux is the same for both
schemes. This is probably due to the inconsistency between the
treatment of inner and boundary cells with these flux limiters. The
other schemes are much more accurate, and very close, even if the
discontinuous Galerkin scheme is a bit more accurate than the
finite volume scheme with slope limiters.}

\section{Conclusions and perspectives}
\label{sec:conclusion}

We have studied two different second order finite volume schemes
applied to a kinetic equation (the BGK model of kinetic gas
theory). We have shown that a scheme with a nonlinear symetric flux limiter is
not second order accurate up to a solid boundary in case of diffuse
reflection conditions, which strongly decreases its accuracy. At the
contrary, the technique of slope limiters allows to reach second order
accuracy up to the boundary. The scheme with slope limiters has been
applied to a 2D supersonic problem with a curvilinear mesh, on which
it shows a much higher accuracy to compute the heat flux along the
wall (as compared to the scheme with flux limiters), which is of
paramount importance in aerodynamics.

Moreover, the finite volume schemes have been compared to an upwind
Discontinuous Galerkin method (with piecewise linear elements). This
scheme does not require any specific treatment at the solid wall, and
the boundary condition can be directly applied to compute the flux at
the gas/solid interface. The
accuracy of this scheme is the same as the finite volume scheme with
slope limiters.
However, the Discontinuous
Galerkin scheme requires twice as degrees of freedom as the finite
volume scheme. Further work would be necessary to compare the
performance of these schemes on 2D curvilinear meshes.

 \bibliographystyle{plain}

 \bibliography{bibliographie}

\newpage

\clearpage

\clearpage

\clearpage

\begin{table}[!!h]
\centering
\begin{tabular}{|c|c|}
\hline
Plates temperature  & 273 K \\ \hline
Left plate velocity  & 0 $\textrm{m.s}^{-1}$ \\ \hline
Right plate velocity  & 300 $\textrm{m.s}^{-1}$ \\ \hline
Distance between the plates  & 1 m \\ \hline
\end{tabular}
\ \ \ \ \ \ 
\begin{tabular}{|c|c|}
\hline
Gas nature & Argon \\ \hline
 $m$ & $0,663 . 10^{-25}$ kg \\ \hline
$R$ & $208,24$ J.kg\textsuperscript{-1}.K\textsuperscript{-1} \\ \hline
$\mu_0$ & $2,117 . 10^{-5}$ Pa.s \\ \hline
$T_0$ & $273,15$ K \\ \hline
$\omega$ & 0,81 \\ \hline
$\alpha$ & 1 \\ \hline
$\textrm{Kn}$ & $9,25 . 10^{-3}$ \\ \hline
\end{tabular}
\caption{{1D Couette flow: flow parameters and gas properties}}
\label{couette}
\end{table}

\clearpage

\clearpage

\clearpage

\clearpage

\clearpage

\clearpage

\end{document}